\documentclass[12pt]{amsart}
\usepackage{amsmath,amssymb,latexsym,amsthm, amscd}

\theoremstyle{definition}
\theoremstyle{plain}
\newtheorem{prop}[subsection]{Proposition}
\newtheorem{thm}[subsection]{Theorem}
\newtheorem{lem}[subsection]{Lemma}
\newtheorem{cor}[subsection]{Corollary}

\newcommand{\mbf}{\mathbf}
\newcommand{\mrm}{\mathrm}

\title[Affine Quivers]{Affine quivers of type $ \tilde{\mbf{A}}_{n}$ and canonical bases}
\author{ Yiqiang Li}
\address{Department of Mathematics\\ Kansas State University\\ Manhattan, KS  66506}
\thanks{}

\begin{document}

\begin{abstract}
Let $Q$ be an affine quiver of type $\tilde{\mbf{A}}_n$. Let $C=(\mbf{c_{ij}})$ 
be the associated generalized Cartan matrix. Let $\mbf{U}^-$ be the negative part of 
the quantized enveloping algebra attached to $C$. In terms of perverse sheaves 
on the moduli space of representations of a quiver, Lusztig constructed 
$\mbf{U}^{\le 0}=\mbf{U}^0 \otimes \mbf{U}^-$ geometrically and gave a canonical basis $B$ of 
$\mbf{U}^-$ at the same time. The simple perverse sheaves which enter $B$ are defined 
abstractly in general. In Lusztig's paper [L2], by using McKay's correspondence, he
gave a description of these canonical basis elements in affine cases, i.e, 
specifying the corresponding supports and local systems. But the chosen 
orientations and the number of vertices of the quiver are not arbitrary. In 
this paper we generalized the description in [L2] for arbitrary orientations 
and vertices  in the case of type $\tilde{\mbf{A}}_n$ by using the theory of 
representations of quivers.
\end{abstract}

\maketitle

\section{Introduction}
Let $C=(c_{\mbf{i},\mbf{j}})_{\mbf{i,j}\in \mbf{I}}$ be a symmetric
generalized Cartan matrix, where $\mbf{I}$ is a finite index set. Let
$\mbf{U}^-$ be the negative part of the quantized enveloping algebra
associated to $C$. A graph $\Gamma$ is attached to $C$ as follows. The
vertex is $\mbf{I}$ and the number of the edges between $\mbf{i}$ and
$\mbf{j}$ ($\mbf{i}\neq\mbf{j}$) is $|c_{\mbf{i,j}}|$ and there is no
edges between $\mbf{i}$ and $\mbf{i}$, for any $\mbf{i}\in\mbf{I}$. 
By indicating the orientation of each
edge in $\Gamma$, we have an oriented graph $\mbf{Q}$, called the
quiver associated to the Cartan matrix $C$.

In [L1], Lusztig studied certain type of simple perverse sheaves on
the representation varieties of $\mbf{Q}$. He constructed
$\mbf{U}^{\leq0}=\mbf{U}^0\otimes\mbf{U}^-$ purely geometrically. In
his construction, the isoclasses of these simple perverse sheaves
produce a canonical basis of $\mbf{U}^-$, with very remarkable
properties. But these simple perverse sheaves are defined
abstractly. We only know that they correspond to orbits in the
representation varieties in the finite case. In [L2], Lusztig
described these simple perverse sheaves in the case of affine type
($\tilde{A}_{n}$, $\tilde{D}_n$, $\tilde{E}_6$, $\tilde{E}_7$ and 
$\tilde{E}_8$). However, his description has some restrictions on the
orientation of the quiver and the number of vertices of the underlying
graph. The orientation he considered are those such that each vertex
is either a sink or a source. 
In the case of type $\tilde{A}_n$
the number of vertices has to be even . 
The theory of representation of affine quivers 
plays a very important role in his description. 
Instead of using the language in [DR] on the
classification of representations of affine quivers, 
Lusztig reinterpreted the problem and redeveloped the theory, 
based on the McKay's correspondence. 

In this paper, we generalize Lusztig's description of the simple
perverse sheaves that enter the canonical basis in the affine case of
type $\tilde{A}_n$. We use the language used in [DR] on the theory of
representations of affine quivers. Note that the language in [DR] does
not have any restriction on the vertices and the orientations, except that
the given orientation cannot be cyclic in the case of type
$\tilde{A}_n$. (The cyclic quiver case has been done by Lusztig.) We
follow the line of what Lusztig did in [L2]. The difficult part is to
prove Lemma 7.5. This leads us to construct a categorical equivalence
between the full subcategory $Rep(T)$, whose objects come from a tube
$T$ of period $p$ and the category of nilpotent representations of a
cyclic quiver with $p$ vertices. The equivalence between these two
abelian categories is well-known. But
our construction gives rise to an isomorphism of certain varieties,
with desired properties that we can use in the proof of Lemma 7.5. 

Section 2 recall the theory of representations of quivers; In Section
3, we study tubes. Section 4, 5 and 6 recall results from
[L2]. Section 7 proves the main theorem. 

I would like to thank my advisor Prof. Zongzhu Lin for his continuing
support and guidance. Without him, this paper will never be
finished. Also I would like to thank Prof. Lusztig and Prof. Ringel for
helpful communications. 

\section{Representations of quivers}

Once and for all, $k$ is an algebraically closed field, of characteristic not necessarily 0.

\subsection{}
A $graph$ $\Gamma$ consists of a pair $(\mbf{I},H)$ with three
operations  
$':H\rightarrow\,\mbf{I}$, $'':H  \rightarrow\,\mbf{I}$ and  $-:H \rightarrow\,H$ satisfying 
the following conditions:
\begin{enumerate}
\item $\mbf{I}$ and $H$ are finite sets;
\item $-:H\rightarrow\,H$ is a fixed point free involution;
\item $\left( \bar{h} \right)'=h'' \quad \text{ and } \quad h' \ne h'' $, for any $h \in H$.
\end{enumerate}

(By abuse of the notation, we write $\Gamma=( \mbf{I}, H)$ for a graph.)
An $orientation$ of a graph $\Gamma=(\mbf{I},H)$ is a subset $\varOmega$ of $H$ such that
\[ \varOmega \cup \bar{\varOmega}=H \text{ and } \varOmega \cap \bar{\varOmega}=\emptyset 
\hspace{.1 cm},\]
where $\bar{\Omega}=\{\bar{h} | h \in \Omega \}$.
A graph $\Gamma= (\mbf{I},H)$ with an orientation $\varOmega$ is called a $quiver$, 
denoted by $\mbf{Q}$=($\Gamma$, $\varOmega$ ). Elements in $\Omega$
are called $arrows$.
Given an arrow $h\in \Omega$, 
we call $h'$ and $h''$ the $starting$ and $terminating$ vertices of $h$, respectively. 
Pictorially, we put $h' \overset{h}{\to} h''$.
A representation of the quiver $\mbf{Q}$ is a pair $(\mbf{V},x)$, 
where $\mbf{V}$ is a finite dimensional $\mbf{I}$-graded $k$-vector space and
$x=\lbrace x_h:\mbf{V}_{h'} \to \mbf{V}_{h''} \hspace{.2 cm}| h \in \varOmega \rbrace$
is a collection of linear maps. 
A morphism $\mbf{\phi}$ between two representations $(\mbf{V},x)$ and
$(\mbf{V}',x')$ is a collection of linear maps
$\mbf{\phi}= \lbrace \phi_{\mbf{i}}:\mbf{V}_{\mbf{i}} \to \mbf{V'}_{\mbf{i}}
\; | \; \mbf{i} \in \mbf{I} \rbrace$
such that $\phi_{h''} x_{h}=x'_{h} \phi_{h'}$ for any $h \in \Omega$. 
Thus, we have just defined an abelian category $Rep(\mbf{Q})$, 
whose objects are representations
of $\mbf{Q}$ and whose morphisms are morphisms of representations of $\mbf{Q}$.

\subsection{}
Let $\mbf{V}$ be an $\mbf{I}$-graded $k$-vector space. 
Let $\mbf{G_V}=\prod_{\mbf{i\in I}} \mrm{GL}(\mbf{V_i})$, 
where $\mrm{GL}(\mbf{V_i})$ is the general linear group of
$\mbf{V_i}$. 
For any subset $H_1 \subseteq H$, define 
\[
\mbf{E}_{\mbf{V}, H_1}=\bigoplus_{h \in H_1}\mrm{Hom}(\mbf{V}_{h'}, \mbf{V}_{h''}).
\]
$\mbf{G_V}$ acts on $\mbf{E}_{\mbf{V}, H_1}$ in the following way. 
For any $g=\{g_{\mbf{i}}|\;\mbf{i\in I} \}$ in $\mbf{G_V}$ and $x=\{x_h|\;h\in H_1\}$ 
in $\mbf{E}_{\mbf{V}, H_1}$, $g.x$ equals $y=\{y_h|\;h \in H_1\}$, 
where $y_h=g_{h''}\circ x_h \circ g_{h'}^{-1}$ for all $h\in H_1$.
In particular, given an orientation $\Omega$ of the graph $\Gamma$, 
\[
\mbf{E}_{\mbf{V},\Omega}=\bigoplus_{h\in \Omega}\mrm{Hom}(\mbf{V}_{h'},\mbf{V}_{h''})
\hspace{.1cm}
\]
and $\mbf{G_V}$ acts on $\mbf{E}_{\mbf{V},\Omega}$. 
Note that for any $x\in \mbf{E}_{\mbf{V},\Omega}$, 
the pair $(\mbf{V},x)$ is a representation of $\mbf{Q}$. 
$(\mbf{V},x)$ and $(\mbf{V},y)$ are isomorphic to each other 
if and only if $x$ and $y$ are in the same $\mbf{G_V}$-orbit. 
Finally, let $\mbf{E_V}=\mbf{E}_{\mbf{V},H}$. 
We call that $x \in \mbf{E_V}$ is $nilpotent$ 
if there exists an $N \geq 2$ such that the following condition is satisfied: 
for any sequence $h_1,\cdots,h_N$ in $H$ such that 
$h'_1=h''_2,\cdots,h'_{N-1}=h''_{N}$, we have $x_{h_1} \circ \cdots \circ x_{h_N}=0$.

\subsection{}
Let $\mbf{Q}=(\Gamma, \Omega)$ be a quiver without oriented cycles. 
A vertex $\mbf{i}$ is called a $sink$ (resp. a $source$) if for any $h\in \Omega$ such that 
$\mbf{i}\in \{h',h''\}$, then $h''=\mbf{i}$ (resp. $h'=\mbf{i}$).
In [BGP], the $Coxeter$ $functors$ $\Phi^+$ and $\Phi^-$ are defined as
follows.  

First, for any $\mbf{i}$ in $\mbf{I}$,
let $\sigma_{\mbf{i}}\Omega= \{h \in \Omega| \mbf{i}\not\in \{h',h''\}\}
\cup \{h\in \bar{\Omega} | \mbf{i}\in \{h',h''\}\}$. 
Then $\sigma_{\mbf{i}}\Omega$ is an orientation of $\Gamma$, 
(i.e, by reversing the arrows in $\Omega$ starting or terminating at
$\mbf{i}$).  
Denote by $\sigma_{\mbf{i}}\mbf{Q}$ the quiver $(\Gamma,\sigma_{\mbf{i}}\Omega)$.

Second, if $\mbf{i}$ is a sink in $\mbf{Q}$, 
define a functor $\Phi_{\mbf{i}}^+:Rep(\mbf{Q}) \to Rep(\sigma_{\mbf{i}}\mbf{Q})$ 
in the following. 
For any $(\mbf{V}, x) \in Rep(\mbf{Q})$,
$\Phi_{\mbf{i}}^+(\mbf{V},x)=(\mbf{W},y)$ is a representation of 
$\sigma_{\mbf{i}}\mbf{Q}$, where $\mbf{W_j}=\mbf{V_j}$
if $\mbf{j}\neq \mbf{i}$ and $\mbf{W_i}$ is the kernel of the linear map
$\sum_{\rho\in \Omega:\;\rho''=\mbf{i}} x_{\rho}\; : 
\oplus_{\rho\in\Omega:\rho''=\mbf{i}} \mbf{V}_{\rho'} \to \mbf{V_i}$; 
for any $h \in \sigma_{\mbf{i}}\Omega$, $y_h=x_h$ if $h'\neq \mbf{i}$,
otherwise if $h'=\mbf{i}$, $y_h:\mbf{W_i}\to \mbf{W}_{h''}$ 
is the composition of the natural maps: 
$\mbf{W_i} \to \oplus_{\rho\in\Omega:\;\rho''=\mbf{i}}\mbf{V}_{\rho'}
=\oplus_{h\in \sigma_{\mbf{i}}\Omega:h'=\mbf{i}}\mbf{W}_{h''}
\to \mbf{W}_{h''}$. 
This extends to a functor. 

Finally, since $\mbf{Q}$ has no oriented cycles, 
we can order the vertices in $\mbf{I}$, 
say $(\mbf{i}_1 ,\cdots,\mbf{i}_n)$ ($|\mbf{I}|=n$), 
in such a way that $\mbf{i}_r$ is a sink
in the quiver $\sigma_{i_{r-1}}\cdots \sigma_{i_1}\mbf{Q}$. 
We then define the Coxeter functor $\Phi^+$ to be 
$\Phi_{\mbf{i}_n}^+\circ \cdots \circ \Phi_{\mbf{i}_1}^+$. 
 
Similarly, if i is a source in $\mbf{Q}$,  
define a functor 
$\Phi_{\mbf{i}}^- :Rep(\mbf{Q}) \to Rep(\sigma_{\mbf{i}}\mbf{Q})$ 
in the following. 
For any $(\mbf{V}, x) \in Rep(\mbf{Q})$, 
$\Phi_{\mbf{i}}^-(\mbf{V}, x)=(\mbf{W}, y)$ is a representation of 
$\sigma_{\mbf{i}}\mbf{Q}$, 
where $\mbf{W_j}=\mbf{V_j}$ if $\mbf{j}\neq \mbf{i}$ and 
$\mbf{W_i}$ equals the cokernel of the linear map 
$\sum_{\rho\in \Omega: \rho'=\mbf{i}}x_{\rho}:
\mbf{V_i} \to \oplus_{\rho\in \Omega: \rho'=\mbf{i}}\mbf{V}_{\rho''}$; 
for any $h \in \sigma_{\mbf{i}}\Omega$, 
$y_{h}=x_{h}$ if $h''\neq \mbf{i}$, 
otherwise if $h''=\mbf{i}$, $y_{h}$ is the composition of the natural map 
$\mbf{W}_{h'} \to \oplus_{\rho\in
  \Omega:\rho'=\mbf{i}}\mbf{V}_{\rho''} \to \mbf{W_i}$.
This extends to a functor. 
The Coxeter functor $\Phi^-$ is defined to be
$\Phi_{\mbf{i}_1}^- \circ \cdots \circ \Phi_{\mbf{i}_{n}}^-$. 

Let $A$ be the path algebra of $\mbf{Q}$. 
We know that  the category, $Rep(\mbf{Q})$, 
of representations of the quiver $\mbf{Q}$ 
is equivalent to the category, $A$-mod, of left $A$-modules 
([CB, \S1 Lemma]) and that the Coxeter functors and the Auslander-Reiten translates are 
equivalent to each other too ([BB]). 
Therefore, the results given in [CB] in the setting of $A$-mod 
can be reinterpreted in the setting of $Rep(\mbf{Q})$. 
We will use some of the results in [CB] directly in the setting of $Rep(\mbf{Q})$.
 
\subsection{}
From now on in this subsection, let $\mbf{Q}=(\Gamma,\Omega)$ 
be an affine quiver without oriented cycles. 

Given an indecomposable representation $(\mbf{V},x)$ of $\mbf{Q}$, 
we call that $(\mbf{V},x)$ is $preprojective$ if $(\Phi^+)^r(\mbf{V},x)=0$ for $r \gg 0$; 
that $(\mbf{V}, x)$ is $preinjective$ if $(\Phi^-)^r(\mbf{V},x)=0$ for $r \gg 0$ 
and that $(\mbf{V},x)$ is $regular$ if $(\Phi^+)^r(\mbf{V},x)\neq 0$ 
and $(\Phi^-)^r(\mbf{V},x)\neq 0$ for any $r>0$. 
We also call that a representation (not necessarily indecomposable) 
is preprojective (resp. preinjective, regular) 
if all indecomposable summands of the representation 
are preprojective (resp. preinjective, regular). 

Let $S_{\mbf{i}}$ be the simple representation corresponding to the vertex $\mbf{i}$. 
It's a representation $(\mbf{V}, x)$ such that $\mbf{V_i}=k$, $\mbf{V_j}=0$ 
if $\mbf{j}\neq \mbf{i}$ and all linear maps $x_h$ are 0. 
(Note that given a graph, the definition of the simple representation 
works for any orientation of the graph. 
By abuse of notation, we always denote by $S_{\mbf{i}}$ 
the simple representation corresponding to vertex $\mbf{i}$ regardless of the orientation.) 

Fix a sequence $(\mbf{i}_1,\cdots,\mbf{i}_n)$ such that $\mbf{i}_r$ is a sink of the quiver
$\sigma_{\mbf{i}_{r-1}}\cdots\sigma_{\mbf{i}_1}(\mbf{Q})$. 
Let $P(\mbf{i}_r)=\Phi_{\mbf{i}_1}^-\circ\cdots\circ\Phi_{\mbf{i}_{r-1}}^-(S_{\mbf{i}_r})$ and
$I(\mbf{i}_r)=\Phi_{\mbf{i}_n}^+\circ\cdots\circ\Phi_{\mbf{i}_{r+1}}^+(S_{\mbf{i}_r})$. 
Note that the set $\{P(\mbf{i}_r)|\;r=1,\cdots,n\}$ 
(resp. $ \{I(\mbf{i}_r)|\;r=1,\cdots,n\}$) gives a complete list of 
projective (resp. injective) representations of $\mbf{Q}$.
Suppose that $(\mbf{V},x)$ is indecomposable, 
then $(\mbf{V},x)$ is preprojective if and only if $(\mbf{V},x )=(\Phi^-)^rP(\mbf{i})$ 
for some $\mbf{i}$ and $r$; 
$(\mbf{V},x)$ preinjective if and only if 
$(\mbf{V},x)=(\Phi^+)^rI(\mbf{i})$ for some $\mbf{i}$ and $r$. 
Furthermore, the regular representations form an extension-closed 
abelian subcategory of $Rep(\mbf{Q})$. 
$\Phi^+$ is an autoequivalence on this subcategory. $\Phi^-$ is its inverse. 
The simple objects in this subcategory are called $regular$ $simples$. 
For each regular simple $(\mbf{V},x)$, there exists a positive integer $r$ such that 
$(\Phi^+)^r(\mbf{V},x)=(\mbf{V},x)$. 
We call the smallest one, $p$, the $period$ of $(\mbf{V},x)$ under $\Phi^+$. 
The set $\{(\mbf{V},x),\cdots,(\Phi^+)^{p-1}(\mbf{V},x)\}$ 
is called the $\Phi^+$-$orbit$ of $(\mbf{V},x)$.
Given a $\Phi^+$-orbit of regular simples, 
the corresponding $tube$, say $T$, 
is a set of all indecomposable regular representations whose 
regular composition factors belong to this orbit. 
We have the following properties:
\begin{itemize}
 \item Every regular indecomposable belongs to a unique tube;
 \item Every indecomposable in a tube has the same period under $\Phi^+$;
 \item All but finitely many regular simples have period one;
 \item The number of non-isomorphic indecomposable representations
       $(\mbf{V},x)$ in  a tube $T$ 
       such that $\underline{\mrm{dim}}\;(\mbf{V},x)=r\delta$ is equal
       to the period of $T$, for any integer $r \geq 1$, where
       $\underline{\mrm{dim}}\;(\mbf{V},x):=\sum_{\mbf{i\in I}}\mrm{dim}\mbf{V_i}\:\mbf{i}$
       is the dimension vector of $(\mbf{V},x)$ and $\delta$ is the
       minimal positive imaginary root of the root system of the underlying graph $\Gamma$.
\end{itemize}
Furthermore we have:
\begin{lem}
If $(\mbf{V},x)$ is regular and  
$\mbf{W} \subseteq \mbf{V}$ is an $\mbf{I}$-graded subspace of $\mbf{V}$ 
such that $\mbf{W}$ is $x$-stable and 
$\underline{\mrm{dim}} (\mbf{V/W})=\underline{\mrm{dim}} \;(\Phi^+)^r R$ 
for some $r$ and $R$ regular indecomposable, then $(\mbf{W},x)$ is regular too. 
(We call that $\mbf{W}$ is $x$-$stable$ if
$x_h(\mbf{W}_{h'})\subseteq\mbf{W}_{h''}$ for any $h\in \Omega$. 
$(\mbf{W},x)$ is the representation induced from
$(\mbf{V},x)$ by restriction.)
\end{lem}
\begin{proof}
Since $\mbf{W}$ is $x$-stable, 
it induces a representation on $\mbf{V/W}$, denote by $(\mbf{V/W},x)$. 
By [DR, Prop. 1.9, 2.6], if $X$ is indecomposable, 
then $X$ is preprojective, regular or preinjective 
if the defect of $X$, $\partial(X)$, is negative, zero or positive. 
We have 
$\partial((\mbf{V},x))=\partial((\mbf{W},x))+\partial((\mbf{V/W},x))$.
Since
$\underline{\mrm{dim}}\;(\mbf{V/W},x)=\underline{\mrm{dim}}\;(\Phi^+)^l(R)$, 
$\partial((\mbf{V/W},x))=0$. 
Thus, $\partial((\mbf{W},x))=0$. 
By [DR, Lemma 3.1], 
we know that every subrepresentation of $(\mbf{V},x)$ have defect $\leq 0$. 
So is every indecomposable summand of $(\mbf{W},x)$. 
Now $(\mbf{W},x)$ can be decomposed into a direct sum of the indecomposables. 
Its defect is then the sum of the defects of these indecomposables.  Hence, 
the defect of every indecomposable summand of $(\mbf{W},x)$ is $0$. Therefore, 
$(\mbf{W},x)$ is regular.
\end{proof}
The following well-known lemma is about the vanishing of $\mrm{Hom}$-groups and
$\mrm{Ext}^1$-groups. Its proof is more or less in [CB].
\begin{lem}
Let $(\mbf{V},x)$ and $(\mbf{W},y)$ be two indecomposable representations of $\mbf{Q}$ 
as follows.
\begin{description}
 \item[a] $(\mbf{V},x)$ and $(\mbf{W},y)$ are both preinjectives such that 
          $(\mbf{V},x)=(\Phi^+)^r I(\mbf{i})$, $(\mbf{W},y)\;=(\Phi^+)^{r'} I(\mbf{i}')$, 
          for some $r\geq r' \geq 0$ and $\mbf{i},\mbf{i}' \in \mbf{I}$;
 \item[b] $(\mbf{V},x)$ and $(\mbf{W},y)$ are both preprojectives such
          that
          $(\mbf{V},x)=(\Phi^-)^{r}P(\mbf{i})$, $(\mbf{W},y)=(\Phi^-)^{r'}P(\mbf{i'})$, 
          for some
          $0 \geq r \geq r'$ and $\mbf{i},\mbf{i}' \in \mbf{I}$;
 \item[c] $(\mbf{V},x)$ and $(\mbf{W},y)$ are both regulars, but they
          are not in the same tube;
 \item[d] $(\mbf{V},x)$ is nonpreinjective, $(\mbf{W},y)$ is preinjective;
 \item[e] $(\mbf{V},x)$ is preprojective, $(\mbf{W},y)$ is nonpreprojective;
\end{description}
then we have
\begin{description}
 \item[A] $\mrm{Ext}^1((\mbf{V},x),(\mbf{W},y))=0$;
 \item[B] $\mrm{Hom}((\mbf{W},y),(\mbf{V},x))=0$; if they are not isomorphic.
\end{description}
\end{lem}

Given two representations $(\mbf{V},x)$ and $(\mbf{W},y)$ of $\mbf{Q}$, 
we have the following exact sequence:
\[
\hspace{-2.5cm}
(1)
\hspace{.5cm}
 0 \to \mrm{Hom}_{\mbf{Q}}((\mbf{V},x),(\mbf{W},y)) \overset{i}{\to}
 \oplus_{\mbf{i \in I}}\mrm{Hom}(\mbf{V_i},\mbf{W_i}) 
 \overset{a}{\to} \oplus_{h \in \Omega}\mrm{Hom}(\mbf{V}_{h'},\mbf{W}_{h''}) 
\]
\[ 
\hspace{11cm}
\overset{b}{\to} \mrm{Ext}^1((\mbf{V},x),(\mbf{W},y)) \to 0,
\]
where $i$ is the inclusion, $a$ and $b$ are defined as follows. Given 
$\theta \in \oplus_{\mbf{i}}\mrm{Hom}(\mbf{V_i},\mbf{W_i})$, 
$a(\theta)=(a(\theta)_h:\mbf{V}_{h'} \to \mbf{W}_{h''}|\; h \in \Omega)$, where 
$a(\theta)_h:=y_h\circ \theta_{h'}-\theta_{h''}\circ x_h$ for any $h\in \Omega$; 
Given $\gamma \in \oplus_h\mrm{Hom}(\mbf{V}_{h'},\mbf{W}_{h''})$, 
$b(\gamma)$ is the following exact sequence:
\[
0 \to (\mbf{W},y) \overset{c}{\to} (\mbf{W}\oplus \mbf{V},
\left[ \begin{array}{cc}
y & \gamma\\
0 & x
\end{array}
\right] ) 
\overset{d}{\to}
(\mbf{V},x) \to 0,
\]
where $c$ and $d$ are inclusion and projection, respectively. (See [R1] for a proof.)
Finally, we have the following Kac's theorem 
\begin{thm}
 \begin{enumerate}
 \item If $(\mbf{V},x)$ is indecomposable then $\underline{\mrm{dim}}\;(\mbf{V},x)$ is a root;
 \item If $\alpha$ is a positive real root, there is a unique (up to isomorphism)
       indecomposable $(\mbf{V},x)$ with 
       $\underline{\mrm{dim}}\;(\mbf{V},x)=\alpha$.
 \end{enumerate}
\end{thm}

\section{Tubes}
From now on, throughout the rest of this paper, we assume that
$\mbf{Q}=(\Gamma, \Omega)$ is an affine quiver of type 
$\tilde{\mbf{A}}_n$ without oriented cycle, unless state explicitly. 

\subsection{}
We know that the minimal positive imaginary root for $\mbf{Q}$ 
is $\delta = \sum_{\mbf{i\in I}}\mbf{i}$. 
Suppose that $T$ is a tube of period $p \neq 1$,
(There are at most 2 such tubes.) 
Fix a regular simple $R$ of $T$ in this section. 
We have $(\Phi^+)^rR\simeq (\Phi^+)^mR$ if $r\equiv m$ (mod $p$).
So it makes sense when we say that the power $r$ is in $\mathbb{Z}/p\mathbb{Z}$.
Denote by $[R]$ the isoclass of the representation R.
Then the set $\{[R], \cdots, [(\Phi^+)^{p-1}R]\}$ is a complete list of 
isoclasses of regular simples in $T$.  
By [CB, Lemma 9.3], we have
$\underline{\mrm{dim}}\; R + \underline{\mrm{dim}}\; \Phi^+ R + 
\cdots +\underline{\mrm{dim}} \;(\Phi^+)^{p-1}R=\delta$.
So $(\underline{\mrm{dim}}\; (\Phi^+)^r R)_\mbf{i}$ equals 0 or 1 
for $\mbf{i} \in \mbf{I}$ and any $r\in \mathbb{Z}/p\mathbb{Z}$. 
Denote by $supp((\Phi^+)^r R)$ the support of the dimension vector of $(\Phi^+)^r R$, 
for any $r\in\mathbb{Z}/p\mathbb{Z}$. 
Thus, we know that the supports of the regular simples in $T$ are 
disjoint. 
For each $r\in\mathbb{Z}/p\mathbb{Z}$, we define a representation $R^T_r$ as follows.
\[ 
\mbf{V_i}=
    \left\{
      \begin{array}{ll}
           k  & \mbox{if $\mbf{i} \in supp((\Phi^+)^r R)$}, \\
           0  & \mbox{otherwise},
      \end{array}
    \right.\; \mbox{for any} \;\mbf{i}\in\mbf{I};
\] 
and
\[ x_h=\left\{\begin{array}{ll}
         1    &\mbox{if $\lbrace h',h'' \rbrace \subseteq supp((\Phi^+)^r R)$},\\
         0    &\mbox{otherwise}, 
              \end{array}
              \right.\;\mbox{ for any} \;h\in \Omega,
\]
where $1$ stands for the identity map.
It's isomorphic to $(\Phi^+)^r R$ and the set
$\{R^T_r|\;r\in\mathbb{Z}/p\mathbb{Z}\}$
gives a complete list of 
the regular simples in $T$ (up to isomorphism). 
Denote by $\mbf{Q}(r)$ the full subquiver of $\mbf{Q}$ 
whose vertex set is $supp((\Phi^+)^r R)$, $r\in \mathbb{Z}/p\mathbb{Z}$. 
Since $(\Phi^+)^rR$ is indecomposable, $\mbf{Q}(r)$ must be connected. So 
$\mbf{Q}(r)$ is either a single point with no arrows or a quiver of the form pictorially: 
\[ 
\bullet \rightarrow \bullet \rightarrow \cdots \rightarrow \bullet.
\]
Denote by $\mbf{s}(r)$ (resp. $\mbf{t}(r)$) the unique source (resp. sink) of $\mbf{Q}(r)$. 
Denote by $\rho(r)$ the unique path from $\mbf{s}(r)$ to $\mbf{t}(r)$ in $\mbf{Q}(r)$. 
Denote by $h(r)$ the unique arrow in $\Omega$ such that
$h(r)'=\mbf{s}(r)$ and $h(r)''=\mbf{t}(r-1)$. 
Given $x \in \mbf{E}_{\mbf{V},\Omega}$, denote by $x_{\rho(r)}$ the natural
composition of the linear transformations on the path $\rho(r)$, 
(i.e, if $\rho(r)$ is a trivial path, 
then $x_{\rho(r)}$ is the identity map; if 
$\rho(r)=h_1 \cdots h_m:
\bullet\overset{h_m}{\to}\bullet\cdots\bullet\overset{h_1}{\to}\bullet$
is a nontrivial path, then 
$x_{\rho(r)}=x_{h_1} \circ \cdots \circ x_{h_m}$ for $m \geq 1$). 
Denote by $Rep(T)$ the full subcategory of $Rep(\mbf{Q})$, 
objects of which are representations whose direct summands are in $T$.

\begin{lem}
Suppose that $(\mbf{V},x)$ is an indecomposable representation in $Rep(T)$, then
\begin{enumerate}
 \item $\mrm{dim}\;\mbf{V_i}=\mrm{dim}\;\mbf{V_j} \;$  if $\mbf{i,j} \in supp((\Phi^+)^r R)$ 
       for some $r$;
 \item $x_h$ is an isomorphism if 
       $\lbrace h',h''\rbrace \subseteq supp((\Phi^+)^r R)$ for some $r$.
\end{enumerate}
\begin{proof}
Note that there exists a sequence of regular representations in $Rep(T)$: 
\[\hspace{-4cm}\star\hspace{4cm} 0=R_0 \subseteq R_1 \subseteq ... \subseteq R_m=(\mbf{V},x),\]
such that $R_1,\;R_2/R_1,\cdots,R_m/R_{m-1}$ are regular simples in $T$ (See [CB]). 
Since for each regular simple, (1) is true, therefore (1) is true for $(\mbf{V},x)$. 

We prove statement (2) by induction on the global dimension of $(\mbf{V},x)$. 
Note that the global dimension of $(\mbf{V},x)$ is defined to be 
$\mrm{dim}\;(\mbf{V},x)=\sum_{\mbf{i\in I}}\mrm{dim}\mbf{V_i}$.
From the construction of the regular simples in $T$,  
the statement is true for any regular simples in $T$. 
Suppose that $(\mbf{V},x)$ is not a regular simple and 
the statement is true for any indecomposable in $Rep(T)$ 
whose global dimension is less than the global dimension of $(\mbf{V},x)$. 
There exists an exact sequence 
\[ 0 \rightarrow (\mbf{V'}, x') \rightarrow (\mbf{V},x) \rightarrow
(\mbf{V''},x'') \rightarrow 0, \]
such that $(\mbf{V'} ,x')$, $(\mbf{V''},x'') \in T$. 
Now for any $h \in \Omega$ such that
$\lbrace h' ,h''\rbrace \subseteq supp ((\Phi^+)^r R)$ for some $r$, 
we have the following commutative diagram:
\[
 \begin{CD}
  0 @>>> \mbf{V'}_{ h' } @>>> \mbf{V}_{ h' } @>>> \mbf{V''}_{h' } @>>> 0\\
  @VVV  @Vx'_hVV @Vx_hVV @Vx''_hVV @VVV\\
  0 @>>> \mbf{V'}_{h''} @>>> \mbf{V}_{h''} @>>> \mbf{V''}_{h''} @>>> 0,
 \end{CD}
\] 
where the horizontal sequences are exact. 
By induction $x'_h$ and $x''_h$ are isomorphisms, 
so by Five lemma, $x_h$ is an isomorphism. (2) holds.
\end{proof}
\end{lem}

\begin{cor}
The Lemma above is true for any representation $(\mbf{V},x) \in Rep(T)$.
\end{cor}
\begin{proof}
Decompose the representation into a direct sum of the indecomposable
representations, apply the Lemma above to each indecomposable summand
of the representation. 
Combine them together, we get the Corollary.
\end{proof}

\begin{prop}
Let $T$ be a tube of period $ p \neq 1$. 
Let $\mbf{Q}_T$ be the cyclic quiver with 
the vertex set $\mathbb{Z}/p\mathbb{Z}$ and the arrow set 
\{$r \rightarrow r-1 |r\in \mathbb{Z}/p\mathbb{Z}$\}. Then 
$Rep(T) \simeq Nil(\mbf{Q}_T)$, 
where $Nil(\mbf{Q}_T)$ is the category of nilpotent representations of $\mbf{Q}_T$. 
\end{prop}
\begin{proof}
Define a functor $F:Rep(T) \rightarrow Nil(\mbf{Q}_T)$ as follows. 
For any $(\mbf{V},x) \in Rep(T)$, 
define a representation $F((\mbf{V},x)):=(F(\mbf{V}),F(x))$ in $Nil(\mbf{Q}_T)$, 
where $F(\mbf{V})=\bigoplus_{r\in \mathbb{Z}/p\mathbb{Z}}F(\mbf{V})_r$
with $F(\mbf{V})_r=\mbf{V}_{\mbf{s}(r)}$ 
for any $r\in \mathbb{Z}/p\mathbb{Z}$ and 
$F(x)=\{F(x)_{r \rightarrow r-1}:F(\mbf{V})_r \to 
F(\mbf{V})_{r-1}|\;r\in\mathbb{Z}/p\mathbb{Z} \}$ 
with 
$F(x)_{r \rightarrow r-1}=x_{\rho(r-1)}^{-1} \circ x_{h(r)}$ 
for $r\in \mathbb{Z}/p\mathbb{Z}$. 
Note that $F(x)$ is well-defined by the Corollary above. 
The nilpotency of $(F(\mbf{V}),F(x))$ comes from 
the existence of a regular sequence $(\star)$  of
$(\mbf{V},x)$ in the proof of Lemma 3.2. This extends to a functor.

Next, we define another functor $G:Nil(\mbf{Q}_T) \rightarrow Rep(T)$ as follows. 
For any $(W,y) \in Nil(\mbf{Q}_T)$, 
define  a representation $G((W,y)):=(G(W),G(y))$ in $Rep(T)$,
where $G(W)=\oplus_{\mbf{i\in I}}G(W)_{\mbf{i}}$ with
$G(W)_{\mbf{i}}=W_{r(\mbf{i})}$ 
if $\mbf{i}\in supp((\phi^+)^{r(\mbf{i})}R$ for some $r(\mbf{i})$ and 
$G(y)=\{G(y)_h:G(W)_{h'}\to G(W)_{h''}|h\in \Omega\}$ with 
$G(y)_h=y_{r\to r-1}$ if $h'=\mbf{s}(r)$, $h''=\mbf{t}(r-1)$ for some $r$,
otherwise $G(y)_h=id$, the identity map.
This extends to a functor too.

Finally, we claim that these two functors are equivalent to each other. 
By definitions, we see that $F \circ G=\mbf{Id}_{Nil(\mbf{Q}_T)}$. 
We construct a natural transformation 
$\alpha:G \circ F \rightarrow \mbf{Id}_{Rep(T)}$ as follows. 
For any $(\mbf{V},x) \in Rep(T)$, $GF((\mbf{V},x))$ is again a
representation in $Rep(T)$ with $GF(V,x)_{\mbf{i}}=\mbf{V}_{\mbf{s}(r)}$, 
where $r=r(\mbf{i})\in\mathbb{Z}/p\mathbb{Z}$ such that 
$\mbf{i}\in supp((\Phi^+)^rR)$, $\mbf{s}(r)$ is the source of the
subquiver $\mbf{Q}(r)$. 
Define a morphism of representations
$\alpha((\mbf{V},x)):GF((\mbf{V},x)) \rightarrow (\mbf{V},x)$ 
by
\[
\alpha((\mbf{V},x))_{\mbf{i}}=x_{\mbf{s}(r) \rightarrow \mbf{i}}:
\mbf{V}_{\mbf{s}(r)}\to \mbf{V_i}
\;\; \mbox{for any} \; \; \mbf{i} \in \mbf{I},
\] 
where $r=r(\mbf{i})$ is such that $\mbf{i} \in supp((\Phi^+)^{r(\mbf{i})} R)$ 
and 
$x_{\mbf{s}(r) \rightarrow \mbf{i}}$ 
is the natural composition of the linear transformations along the path 
$\rho_{\mbf{s}(r) \rightarrow \mbf{i}}$ 
(i.e. the path from $\mbf{s}(r)$ to $\mbf{i}$ in $\mbf{Q}(r)$).
This really is a morphism of representatios of $\mbf{Q}$. 
It extends to a natural transformation.
By Corollary 3.3, we see that $\alpha((\mbf{V},x))_{\mbf{i}}$ is
isomorphic, for any $\mbf{i}\in\mbf{I}$.
Therefore, $\alpha$ is a natural isomorphism.
The Proposition is proved.
\end{proof}

\subsection{}
From the analysis in Section 3.1, we know that the arrow set 
$\Omega$ can be separated into two disjoint subsets, written as 
$\Omega=\Omega_1\cup \Omega_2$, where
$\Omega_1$ consists of all arrows, $h$, 
such that $\{h',h''\}\subseteq supp((\Phi^+)^rR)$ for 
some $r$ and $\Omega_2=\{h(r)|\;r\in \mathbb{Z}/p\mathbb{Z}\}$. 
By definition, 
\[
\mbf{E}_{\mbf{V},\Omega}=\mbf{E}_{\mbf{V},\Omega_1}\oplus \mbf{E}_{\mbf{V},\Omega_2}.
\]
Let $S$ be the subset of $\mbf{E}_{\mbf{V},\Omega}$ 
consisting of all elements $x\in\mbf{E}_{\mbf{V},\Omega}$ such that 
the representation $(\mbf{V},x)$ is in $Rep(T)$. 
By definition, $S$ is $\mbf{G_V}$-invariant. If $S$ is nonempty, 
then by Corollary 3.3, we know that $\mrm{dim}\mbf{V_i}=\mrm{dim}\mbf{V_j}$ if 
$\{\mbf{i,j}\}\subseteq supp((\Phi^+)^rR)$ for some $r$. 
Since $\mbf{E}_{\mbf{V},\Omega}\cong \mbf{E}_{\mbf{V}',\Omega}$ 
if $\mbf{V}\cong\mbf{V}'$, 
without lost of generality, we may assume that $\mbf{V_i}=\mbf{V_j}$ if 
$\{\mbf{i,j}\}\subseteq supp((\Phi^+)^rR)$ for some $r$. 
Denote by $\mrm{Aut}_{\mbf{V},\Omega_1}$ the (open) subset of $\mbf{E}_{\mbf{V},\Omega_1}$ 
consisting of all elements whose components are isomorphisms. 
By Corollary 3.3, $S$ is contained in 
$\mrm{Aut}_{\mbf{V},\Omega_1}\oplus\mbf{E}_{\mbf{V},\Omega_2}$. 
From the construction of the functors
$F$ and $G$ in the proof of Proposition 3.4, we define two maps 
$F:\mrm{Aut}_{\mbf{V},\Omega_1}\oplus\mbf{E}_{\mbf{V},\Omega_2} 
\to \mbf{E}_{F(\mbf{V}),\Omega_T}$ 
and 
$G:\mbf{E}_{F(\mbf{V}), \Omega_T}
\to \mrm{Aut}_{\mbf{V},\Omega_1}\oplus\mbf{E}_{\mbf{V},\Omega_2}$
as follows. 
$F(x)_{r\to r-1}=x_{\rho(r-1)}^{-1}\circ x_{h(r)}$, 
for any $x\in \mrm{Aut}_{\mbf{V},\Omega_1}\oplus\mbf{E}_{\mbf{V},\Omega_2}$,
$r\in \mathbb{Z}/p\mathbb{Z}$; 
$G(y)_h=y_h$ if $h\in \Omega_2$, $G(y)_h=id$, the identity map, if $h\in \Omega_1$, 
for any $y\in \mbf{E}_{F(\mbf{V}),\Omega_T}$, respectively.
Note that the image of $S$ under $F$ is the subset of $\mbf{E}_{F(\mbf{V}),\Omega_T}$ 
consisting of all nilpotent elements. 
Let $S'$ be the image of $F(S)$ under $G$ (i.e, $S'=GF(S)$). 
We have $S'\subseteq S$, since $F(S)$ consists of nilpotent elements. 
Furthermore, by construction, $S'\subseteq \{id\}\oplus \mbf{E}_{\mbf{V},\Omega_2}$, 
where $id$ stands for the element in $\mbf{E}_{\mbf{V},\Omega_1}$
whose components are identity maps.

(a) The restriction 
$F|_{\{id\}\oplus\mbf{E}_{\mbf{V},\Omega_2}}: \{id\}\oplus\mbf{E}_{\mbf{V},\Omega_2} \to
\mbf{E}_{F(\mbf{V}),\Omega_T}$ 
is obviously an isomorphism.

Its inverse is 
$G:\mbf{E}_{F(\mbf{V}),\Omega_T}\to \{id\}\oplus\mbf{E}_{\mbf{V},\Omega_2}$.
Let $\mbf{H_V}$ be the stabilizer of $\{id\}\oplus\mbf{E}_{\mbf{V},\Omega_2}$ in $\mbf{G_V}$.
In fact, 
\[
\mbf{H_V}=\{g=(g_{\mbf{i}})\in \mbf{G_V}|\; g_{\mbf{i}}=g_{\mbf{j}}, \;\mbox{if}\; 
\{\mbf{i,j}\}\subseteq supp((\Phi^+)^rR) \;{\emph for some }\; r\}.
\]   
So we have an action of $\mbf{H_V}$ on
$\{id\}\oplus\mbf{E}_{\mbf{V},\Omega_2}$, induced from the action of
$\mbf{G_V}$ on $\mbf{E}_{\mbf{V},\Omega}$. 
By definition, the projection $F:\mbf{H_V}\to \mbf{G}_{F(\mbf{V})}$ is an isomorphism.
Thus the action is compatible with the action of $\mbf{G}_{F(\mbf{V})}$ on 
$\mbf{E}_{F(\mbf{V}),\Omega_T}$. 
(I.e, $F(gx)=F(g)F(x)$, for any $g\in \mbf{H_V}$ 
and $x\in\{id\}\oplus\mbf{E}_{F(\mbf{V}),\Omega_T}$.) 
Now define an action of $\mbf{H_V}$ on 
$\mbf{G_V}\times(\{id\}\oplus\mbf{E}_{\mbf{V},\Omega_2})$
by $(g,x).g'=(gg',g'^{-1}x)$, 
for any 
$(g,x)\in \mbf{G_V}\times(\{id\}\oplus\mbf{E}_{\mbf{V},\Omega_2})$ 
and $g'\in \mbf{H_V}$.
The action of $\mbf{G_V}$ on $\mbf{E}_{\mbf{V},\Omega}$ defines a map
$m:\mbf{G_V}\times(\{id\}\oplus\mbf{E}_{\mbf{V},\Omega_2})\to 
\mrm{Aut}_{\mbf{V},\Omega_1}\oplus\mbf{E}_{\mbf{V},\Omega_2}$. 
Moreover, we have
\begin{lem}
\begin{enumerate}
\item $(\mbf{G_V}\times(\mrm{Aut}_{\mbf{V},\Omega_1}\oplus\mbf{E}_{\mbf{V},\Omega_2}),m)$ 
      is a quotient of $\mbf{G_V}\times(\{id\}\oplus\mbf{E}_{\mbf{V},\Omega_2})$ by 
      $\mbf{H_V}$. $($See $[B, 6.3]$ for a definition of a quotient of
      a variety by an algebraic group.$)$ We write
      $\mbf{G_V}\times^{\mbf{H_V}}(\{id\}\oplus\mbf{E}_{\mbf{V},\Omega_2})=
      \mrm{Aut}_{\mbf{V},\Omega_1}\oplus\mbf{E}_{\mbf{V},\Omega_2}$.
\item $\mbf{G_V}\times^{\mbf{H_V}} S'= S$ and  
      $\mbf{G_V}\times^{\mbf{H_V}} O'= O_x$,
      where $O_x$ is the $\mbf{G_V}$-orbit of $x$ in 
      $\mrm{Aut}_{\mbf{V},\Omega_1}\oplus\mbf{E}_{\mbf{V},\Omega_2}$ and 
      $O'$ is the $\mbf{H_V}$-orbit of $GF(x)$ in $\{id\}\oplus\mbf{E}_{\mbf{V},\Omega_2}$.
\end{enumerate}
\end{lem}
\begin{proof}
To prove (1), by Proposition 6.6 in [B], it suffices to prove the
following three conditions. 
\begin{itemize}
\item [(a)] $m$ is a separable orbit map;
\item [(b)] $\mrm{Aut}_{\mbf{V},\Omega_1}\oplus\mbf{E}_{\mbf{V},\Omega_2}$ is normal; 
\item [(c)] The irreducible components of 
            $\mbf{G_V}\times(\{id\}\oplus\mbf{E}_{\mbf{V},\Omega_2})$ are open.
\end{itemize}
Since $\mbf{G_V}$ and $\mbf{E}_{\mbf{V},\Omega_2}$ are irreducible, 
so is $\mbf{G_V}\times(\{id\}\oplus\mbf{E}_{\mbf{V},\Omega_2})$. (c) holds. 
Since $\mrm{Aut}_{\mbf{V},\Omega_1}\oplus\mbf{E}_{\mbf{V},\Omega_2}$ is a smooth variety, 
it's normal. (b) holds.
 
Fix $y$ in $\mrm{Aut}_{\mbf{V},\Omega_1}\oplus\mbf{E}_{\mbf{V},\Omega_2}$. 
Given $(g,x)$ and $(g',x')$ in $m^{-1}(y)$,  we have $gx=g'x'$. 
So $g'^{-1}gx=x'$. 
Since $x$ and $x'$ are in $\{id\}\oplus\mbf{E}_{\mbf{V},\Omega_2}$ and 
$\mbf{H_V}$ is the stabilizer of $\{id\}\oplus\mbf{E}_{\mbf{V},\Omega_2}$ in $\mbf{G_V}$,
$g'^{-1}g$ is in $\mbf{H_V}$. 
Let $f=g^{-1}g'$, we have $(g,x)f=(gf, f^{-1}x)=(g',x')$.
So $(g,x)$ and $(g',x')$ are in the same $\mbf{H_V}$-orbit.
By definition, $m((g,x)f)=m((g,x))$, for any 
$(g,x)\in \mbf{G_V}\times(\{id\}\oplus\mbf{E}_{\mbf{V},\Omega_2})$ 
and $f\in \mbf{H_V}$.
So $m^{-1}(y)$ is an $\mbf{H_V}$-orbit. 
On the other hand, given $y$ in
$\mrm{Aut}_{\mbf{V},\Omega_1}\oplus\mbf{E}_{\mbf{V},\Omega_2}$,
the natural transformation $\alpha$ in the proof of Proposition 3.4 gives 
an isomorphism $(\mbf{V},y)\simeq (\mbf{V},G(F(y)))$. 
So there exists $g$ in $\mbf{G_V}$ such that $y=gG(F(y))$. 
Thus, $m((g,G(F(y))))=y$. 
Hence $m$ is surjective. 
Therefore, $m$ is an orbit map. 
To prove (a), it's enough to prove that $m$ is separable. 
Fix a point $(1,x)$ in $\mbf{G_V}\times(\{id\}\oplus\mbf{E}_{\mbf{V},\Omega_2})$, 
where $1$ stands for the element in $\mbf{G_V}$ whose components are identity maps. 
$m$ induces a map of tangent spaces
\[
(\mrm{d} m)_{(1,x)}: gl_{\mbf{V}}\times (\{id\} \oplus \mbf{E}_{\mbf{V},\Omega_2})\to 
\mbf{E}_{\mbf{V},\Omega}.
\]   
It's defined by $(\mrm{d}m)_{(1,x)}(g,y)=z \;(=(z_h)_{h\in \Omega})$, 
where $z_h=y_h+g_{h''}x_h-x_hg_{h'}$ for any $(g,y)$ in 
$\mbf{G_V}\times(\{id\}\oplus\mbf{E}_{\mbf{V},\Omega_2})$.
By [B, Theorem 17.3], to prove that $m$ is separable, it's enough to prove that 
$(\mrm{d}m)_{(1,x)}$ is surjective. 
Given $z\in \mbf{E}_{\mbf{V},\Omega}$, let $y_h=z_h$ if $h\in\Omega_2$ 
and $y_h=id$ if $h\in \Omega_1$. 
$y=(y_h)_{h\in \Omega}$ is in $\{id\}\oplus\mbf{E}_{\mbf{V},\Omega}$, by definition.
Note that the set $\Omega_1$ can be decomposed into disjoint union of subsets of arrows 
such that arrows in each subset can be ordered as follows.
\[\mbf{i}_1\overset{h_1}\to \mbf{i}_2\overset{h_2}\to\cdots\overset{h_n}\to\mbf{i}_{n+1}.\]
Let $\phi_{\mbf{i}_1}=id$, $\phi_{\mbf{i}_2}=z_{h_1}$, 
$\phi_{\mbf{i}_r}=z_{h_{r-1}}+z_{h_{r-2}}+\cdots+z_{h_1}-(r-2)$ for $2\leq r\leq n+1$. 
This defines an element, $\phi$, in $gl_{\mbf{V}}$. We have $(\mrm{d}m)_{(1,x)}(\phi, y)=z$. 
So (a) holds. (1) is proved.

Note that the stabilizer of $S'$ (resp. $O'$) in $\mbf{G_V}$ is also $\mbf{H_V}$, 
$m(\mbf{G_V}\times S')=S$ (resp. $m(\mbf{G_V}\times O')=O_x$). By (1), we have (2).
\end{proof}
\section{preparatory results}
\subsection{}
We recall some definitions and results from [L1] and [L2].
Let $\mbf{Q}=(\Gamma, \Omega)$ be a quiver.
Let $\mbf{V}$ be a finite dimensional 
$\mbf{I}$-graded vector space over $k$. Let 
$\mbf{E_V}=\oplus_{h \in H} \mrm{Hom}_k(\mbf{V}_{h'},\mbf{V}_{h''})$ and 
$\mbf{G_V}=\prod_{\mbf{i \in I}} GL(\mbf{V_i})$, as in Section 2.2.
The Lie algebra of $\mbf{G_V}$ is 
$gl_{\mbf{V}}=\oplus_{\mbf{i}\in \mbf{I}} \mrm{End}(\mbf{V_i})$. 
Recall that $\mbf{G_V}$ acts on $\mbf{E_V}$ by $(g,x)=y$, 
where $y_h=g_{h''}x_h g^{-1}_{h'}$ for any $h \in H$. 
On $\mbf{E_V}$, a non-degenerate symplectic form $<,>$ is defined by 
\[
<x,x'>=\sum_{h \in \Omega} \mrm{tr}(x_h x'_{\bar{h}})-
\sum_{h \in \bar{\Omega}}\mrm{tr}(x_h x'_{\bar{h}}),
\] 
where $x,\;x'\in\mbf{E_V}$ and  $\mrm{tr}$ is the trace of the endomorphisms of 
$\mbf{V}_{h''}$. 
This symplectic form is $\mbf{G_V}$-invariant.
Since $\mbf{E_V}=\mbf{E}_{\mbf{V},\Omega}\oplus \mbf{E}_{\mbf{V},\bar{\Omega}}$,
$\mbf{E}_{\mbf{V},\Omega}$ and $\mbf{E}_{\mbf{V},\bar{\Omega}}$ appear as complementary 
Lagrangian subspaces of $\mbf{E_V}$. 
(They are $\mbf{G_V}$-stable.) 
$<,>$ defines a non-singular pairing 
$\mbf{E}_{\mbf{V},\Omega} \otimes \mbf{E}_{\mbf{V},\bar{\Omega}} \rightarrow k$.
So $\mbf{E_V}$ can be regarded as the cotangent bundle of $\mbf{E}_{\mbf{V},\Omega}$.
In particular, if $Y$ is a subvariety of $\mbf{E}_{\mbf{V},\Omega}$, 
then the conormal bundle of $Y$ may be naturally regarded as a subvariety of $\mbf{E_V}$. 

Attached to the $\mbf{G_V}$-action on the symplectic vector space 
$\mbf{E_V}$, a moment map $\Psi:\mbf{E_V} \rightarrow gl_{\mbf{V}}$, where the 
$\mbf{i}^{th}$-component $\Psi_{\mbf{i}}:\mbf{E_V} \rightarrow \mrm{End}(\mbf{V_i})$ 
is given by
\[
\Psi_{\mbf{i}}(x)=\sum_{h \in \Omega:h''=\mbf{i}}x_h x_{\bar{h}}
-\sum_{h\in \bar{\Omega}:h''=\mbf{i}} x_h x_{\bar{h}}.
\]
We have

(1) If $x'\in \mbf{E}_{\mbf{V},\Omega}$ , $x'' \in \mbf{E}_{\mbf{V},\bar{\Omega}}$ then 
$\Psi(x'+x'')=0$ if and only if $x''$ is orthogonal with respect to
$<,>$ to the tangent space at $x'$
to the $\mbf{G_V}$-orbit of $x'$.

Given $x\in \mbf{E}_{\mbf{V},\Omega}$, 
let $\mbf{T}$ be the tangent space at $x$ to the $\mbf{G_V}$-orbit of $x$ and
let $\mbf{T'}$ be the set of vectors in $\mbf{E}_{\mbf{V},\bar{\Omega}}$ 
which are orthogonal to $\mbf{T}$ under $<,>$.
Consider the exact sequence 2.4 (1) when $(\mbf{V},x)=(\mbf{W},y)$:
\[
\hspace{-2.5cm}
(2)
\hspace{.5cm}
 0 \to \mrm{Hom}((\mbf{V},x),(\mbf{V},x)) \to
 \oplus_{\mbf{i \in I}}\mrm{Hom}(\mbf{V_i},\mbf{V_i}) 
 \overset{a}{\to} \oplus_{h \in \Omega}\mrm{Hom}(\mbf{V}_{h'},\mbf{V}_{h''}) 
\]
\[ 
\hspace{11cm}
\overset{b}{\to} \mrm{Ext}^1((\mbf{V},x),(\mbf{V},x)) \to 0.
\]
By definitions, we have that $\mbf{T}$ is the image of $a$. 
Since the sequence is exact, the image of $a$ is the kernel of $b$, 
so 
\[
\hspace{-6.5cm}
(3) 
\hspace{3cm}
\mbf{T'}\simeq \mrm{Ext}^1((\mbf{V},x),(\mbf{V},x))^*.
\]

\subsection{}
Define $\Lambda_{\mbf{V}}$ to be the set of all nilpotent elements  
$x$ in $\mbf{E_V}$ such that $\Psi(x)=0$. 
By Theorem 12.3 in [L1],
$\Lambda_{\mbf{V}}$ is a closed $\mbf{G_V}$-stable subvariety of
$\mbf{E_V}$ of pure dimension $\mrm{dim}\mbf{(E_V)}/2$. 
(Pure means that all irreducible components have the same dimension.) 
Let $\mbf{Irr}\Lambda_{\mbf{V}}$ be the set of 
all irreducible components of $\Lambda_{\mbf{V}}$. 

If $\mbf{W}$ is a subspace of $\mbf{V}$ and $x \in \mbf{E_V}$, recall that 
$\mbf{W}$ is $x$-$stable$ if $x_h (\mbf{W}_{h'}) \subseteq \mbf{W}_{h''}$ 
for any $h \in H$. 
Let $\mathcal{X}$ be the set of all sequences
$\boldsymbol{\nu}=(\nu^1,...,\nu^m)$, 
where $\nu^r \in \mathbb{N}[\mbf{I}]$ for $r=1,\cdots, m$. 
Let $\mathcal{Y}$ be the subset of $\mathcal{X}$ such that 
$\nu^r$ is discrete, for all $r$.
($\nu^r$ is discrete if  $\nu^r_{h'} \cdot \nu^r_{h''}=0$, for any $h \in H$.)  

Given $\boldsymbol{\nu} \in \mathcal{X}$,
a flag $\mbf{f}$ of type $\boldsymbol{\nu}$ is, by definition, a sequence 
$\mbf{f}=(\mbf{V}=\mbf{V}^0 \supseteq \mbf{V}^1 \supseteq \cdots \supseteq \mbf{V}^m =0)$
of $\mbf{I}$-graded subspaces of $\mbf{V}$ such that 
$\underline{\mrm{dim}} \;\mbf{V}^{r-1}/\mbf{V}^r = \nu^r$ 
for all $r=1,\dots,m-1$. 
If $x \in \mbf{E_V}$, we call that $\mbf{f}$ is $x$-$stable$ 
if $\mbf{V}^r$ is $x$-stable for $r=0,\dots,m-1$.

Given a sequence $\boldsymbol{\nu} \in \mathcal{Y}$ such that 
$|\boldsymbol{\nu}|=\underline{dim}\;\mbf{V}$. 
($|\boldsymbol{\nu}|=\sum_{r=1}^m\nu^r$, if 
$\boldsymbol{\nu}=(\nu^1,\cdots,\nu^m)$.)
Define a function 
$\chi_{\boldsymbol{\nu}}:\Lambda_{\mbf{V}} \rightarrow \mathbb{Z}$ as follows.
For any $x\in\Lambda_{\mbf{V}}$, $\chi_{\boldsymbol{\nu}}(x)$ is 
the Euler characteristic of the variety of all $x$-stable flags 
of type $\boldsymbol{\nu}$ in $\mbf{V}$. 
Let $\emph{F}_{\mbf{V}}$ be the $\mathbb{Q}$-vector space spanned by the functions 
$\chi_{\boldsymbol{\nu}}$ for all $\boldsymbol{\nu} \in \mathcal{Y}$
such that $|\boldsymbol{\nu}|=\underline{\mrm{dim}}\mbf{V}$.
From [L1, \S 12], 
we know that $\mathit{F}_{\mbf{V}}$ is a finite dimensional vector space 
and all functions in $\mathit{F}_{\mbf{V}}$ are constructible.
Furthermore, we have 
\begin{lem}
$|\mbf{Irr} \Lambda_{\mbf{V}}| \leq \mrm{dim}\; \mathit{F}_{\mbf{V}}$,
for any $\mbf{I}$-graded $k$-vector space $\mbf{V}$.
\end{lem}
\begin{proof}
For any $Y \in \mbf{Irr}\Lambda_{\mbf{V}}$, 
we can find an open dense subset $Y_0$ of $Y$ such that 
any function $\chi \in \mathit{F}_{\mbf{V}}$ is constant in $Y_0$. 
This is because any function in $\mathit{F}_{\mbf{V}}$ is constructible and 
hence for each function in $\mathit{F}_{\mbf{V}}$, 
there exists such open dense subset in $Y$. On the other hand, 
$\mathit{F}_{\mbf{V}}$ is finite dimensional, 
hence we can choose a basis $\chi_r$ of $\mathit{F}_{\mbf{V}}$
and for each element $\chi_r$ in the basis, 
choose an open dense subset, say $Y_r$, of $Y$ such that $\chi_r$ is constant in $Y_r$. 
Now let $Y_0=\cap_{r} Y_r$. $Y_0$ satisfies the required condition. 
Let $\mathbb{Q}[\mbf{Irr}\Lambda_{\mbf{V}}]$ be the $\mathbb{Q}$-vector space 
with basis being the elements in $\mbf{Irr}\Lambda_{\mbf{V}}$. 
Define a linear map 
$\theta : \mathit{F}_{\mbf{V}} \to \mathbb{Q}[\mbf{Irr}\Lambda_{\mbf{V}}]$
as follows. 
For any $\chi \in \mathit{F}_{\mbf{V}}$,   
$\theta(\chi):=\sum_{Y\in \mbf{Irr}\Lambda_{\mbf{V}}} \chi(Y_0)\;Y$, 
where $\chi(Y_0)$ stands for the constant value of $\chi$ on $Y_0$. 
But for any $Y \in \mbf{Irr}\Lambda_{\mbf{V}}$, 
by [L2, 3.6; L4, Lemma 2.4], 
there exists a function $\chi_Y \in \mathit{F}_{\mbf{V}}$ with the following properties: (1) 
there exists an open dense subset $Y_0$ in $Y$ such that $\chi_Y|_{Y_0} \equiv 1$ and (2) 
there exists a closed subset $\Delta$ in $\Lambda_{\mbf{V}}$ of dimension strictly less than 
the dimension of $\Lambda_{\mbf{V}}$ such that 
$\chi_Y|_{\Lambda_{\mbf{V}}-Y\cup\Delta} \equiv 0$. 
We then have $\theta(\chi_Y)=Y$. So $\theta$ is surjective. The Lemma is proved.
\end{proof}

\subsection{}
Let $U^-$ be the negative part of the enveloping algebra associated to the
symmetric generalized Cartan matrix $C$ of type $\tilde{A}_n$. 
When $n=1$, $U^-$ is the $\mathbb{Q}$-algebra with two generators
$f_{\mbf{i}}$ and $f_{\mbf{j}}$ and two relations:
\[f_{\mbf{i}}\;f^3_{\mbf{j}}-3f_{\mbf{j}}\;f_{\mbf{i}}\;f^2_{\mbf{j}}
+3f^2_{\mbf{j}}\;f_{\mbf{i}}\;f_{\mbf{j}}-f^3_{\mbf{j}}\;f_{\mbf{i}}=0,
\]
\[
f_{\mbf{j}}\;f^3_{\mbf{i}}-3f_{\mbf{i}}\;f_{\mbf{j}}\;f^2_{\mbf{i}}
+3f^2_{\mbf{i}}\;f_{\mbf{j}}\;f_{\mbf{i}}-f^3_{\mbf{i}}\;f_{\mbf{j}}=0.
\]
When $n\geq 2$, $U^-$ is the $\mathbb{Q}$-algebra with generators
$f_{\mbf{i}}$, where $\mbf{i\in I}$ ($|\mbf{I}|=n+1$), and relations:
\[
f_{\mbf{i}}\;f^2_{\mbf{j}}-2f_{\mbf{i}}\;f_{\mbf{j}}\;f_{\mbf{i}}
+f^2_{\mbf{j}}\;f_{\mbf{i}}=0,
\;\;\emph{if}\; c_{\mbf{i,j}}=-1;
\] 
\[f_{\mbf{i}}\;f_{\mbf{j}}=f_{\mbf{j}}\;f_{\mbf{i}}, 
\;\;\emph{if}\;c_{\mbf{i,j}}=0.\]
 
For any $\nu\in\mathbb{N}[\mbf{I}]$, fix an $\mbf{I}$-graded
$k$-vector space, $\mbf{V}_{\nu}$, of dimension vector $\nu$.
Consider the $\mathbb{Q}$-vector space 
$\emph{F}=\bigoplus_{\nu\in\mathbb{N}[\mbf{I}]}\emph{F}_{\mbf{V}_{\nu}}$.
It's indenpendent of the choice of $\mbf{V}_{\nu}$, 
since $\emph{F}_{\mbf{V}_{\nu}}\simeq \emph{F}_{\mbf{V'}_{\nu}}$,
for any $\mbf{V}_{\nu}$ and $\mbf{V'}_{\nu}$ of the same dimension vector.
By [L1, Theorem 12.13], there is an algebra structure on $\emph{F}$
such that the map 
$(f^{r_1}_{\mbf{i}_1}/r_1!) \cdots (f^{r_m}_{\mbf{i}_m} / r_m!) 
\mapsto \chi_{\boldsymbol{\nu}}$, 
where $\boldsymbol{\nu}=(r_1 \mbf{i}_1,\cdots,r_m \mbf{i}_m)$, defines
an algebra isomorphism $U^-\to \emph{F}$. 
There is a natural grading
$U^-=\oplus_{\nu\in\mathbb{N}[\mbf{I}]}U^-_{\nu}$ and the isomorphism
respects the grading. 
Consequently, we have 
\begin{lem}
$\mrm{dim}\;\mathit{F}_{\mbf{V}} = \mrm{dim} \; U^-_{\nu}$,
for any $\mbf{V}$ of dimension vector $\nu$.
\end{lem}

\section{Preparatory results 2}
\subsection{}
Given  an $\mbf{I}$-graded vector space $\mbf{V}$ over $k$, 
we consider the set $\varphi(\mbf{V},\Omega)$ of all pairs
$(\sigma,\lambda)$, where $\sigma: \mbf{Ind} \to \mathbb{N}$ 
is a function and 
$\lambda$ is $(0)$ or a sequence of decreasing positive integers 
$(\lambda_1 \geq \lambda_2 \geq \cdots \geq \lambda_q)$, 
satisfying the following conditions:
\begin{enumerate}
\item $\sigma$ has finite support;
\item $\prod^{p-1}_{r=0} \sigma([(\Phi^+)^r R])=0$, for any isoclass
      $[R]\in\mbf{Ind}$ where $R$ is a regular indecomposable
      representation of period  $p \geq 2$;
\item $\sigma([R])=0$, for any isoclass $[R]\in\mbf{Ind}$ where $R$ is
      a regular indecomposable representation of period $1$;
\item $\sum_{[P] \in \mbf{Ind}} \sigma([P]) \underline{\mrm{dim}}\; P +
      \sum_r \lambda_r\; \delta = \underline{\mrm{dim}}\;\mbf{V}$.
\end{enumerate}

Given $(\sigma,\lambda) \in \varphi(\mbf{V},\Omega)$, 
we consider  the subset $\mbf{X}(\sigma,\lambda)$ of 
$\mbf{E}_{\mbf{V},\Omega}$ consisting of all elements $x$ such that 
\[
(\mbf{V},x) \simeq \bigoplus_{[P]\in\mbf{Ind}} P^{\sigma([P])} 
\oplus R_1\oplus \cdots \oplus R_q
\hspace{.1cm},
\]
where 
$P^{\sigma([P])}$ is the direct sum of $\sigma([P])$ copies of $P$ and
$R_1,\cdots,R_q$ are regular indecomposables 
from different tubes of period $1$ such that 
$\underline{\mrm{dim}}R_r=\lambda_r\delta$ for $r=1,\cdots,q$. We have

\begin{prop}
 \begin{enumerate}
  \item $\mbf{X}(\sigma,\lambda)$ is open dense smooth in its closure 
        $\overline{\mbf{X}(\sigma,\lambda)}$. It's also irreducible of dimension
        equal to $q+\mrm{dim}\;O_x$, where $O_x$ is the $\mbf{G_V}$-orbit of $x$ 
        in $\mbf{X}(\sigma,\lambda)$;
  \item Let $N(\sigma,\lambda)$ be the conormal bundle of 
        $\overline{\mbf{X}(\sigma,\lambda)} \subseteq \mbf{E}_{\mbf{V},\Omega}$, 
        regarded as a subvariety of $\mbf{E_V}$, then $N(\sigma,\lambda)$ is an
        irreducible component of $\Lambda_{\mbf{V}}$;
  \item For any two different pairs $(\sigma,\lambda)$ and
        $(\sigma',\lambda')$ in $\varphi(\mbf{V},\Omega)$, then 
        $N(\sigma,\lambda) \neq N(\sigma',\lambda')$.   
 \end{enumerate}
\end{prop}  
For a proof of (1) (resp. (2)), see [R2, Theorem 4.3] 
(resp. [R2, Corollary 5.3]). (3) is from (2). 
Also see [L2, Proposition 4.14] for a proof of this Proposition 
with some restriction on the orientation of the quiver.

\subsection{}  
Fix an element $x \in \mbf{X}(\sigma,\lambda)$, 
the representation $(\mbf{V},x)$ can be decomposed into a direct sum of indecomposables.  
Rewrite $\mbf{V}$ as a direct sum $\mbf{V}=\oplus_{r\in \mathbb{Z}}\mbf{V}(r)$, 
where $\mbf{V}(r)$ is an $\mbf{I}$-graded subspace of $\mbf{V}$ for any $r$, 
such that  $\mbf{V}(r)$ is $x$-stable.
There exists some integers $n_1 \leq n_2$ such that we can order
$V(r)$ in the following way.   
\begin{enumerate}
 \item For any $r$ such that $r \leq n_1$, there exists $s=s(r) \geq 0$ and 
       $\mbf{i}=\mbf{i}(r) \in \mbf{I}$ such that $(\mbf{V}(r),x)$ is isomorphic to  a direct 
       sum of copies of $(\Phi^-)^s P(\mbf{i})$, where $P(\mbf{i})$ is the projective 
       representation 
       corresponding to $\mbf{i}$; moreover, if $r<r'\leq n_1$, then either $s(r) < s(r')$ 
       or $s(r)=s(r')$ and $\mbf{i}(r) \neq \mbf{i}(r')$.
 \item For any $r$ such that $n_1<r\leq n_2$, $(\mbf{V}(r),x)$ are regulars in different tubes
       of period $\geq 2$; 
       For $r=n_2+1$, $(\mbf{V}(r),x) \simeq R_1\oplus \cdots \oplus R_q$,
       where $R_1,\cdots,R_q$ are regular indecomposables from different tubes of period $1$ 
       such that $\underline{\mrm{dim}}R_s=\lambda_s\delta$ for
       $s=1,\cdots,q$. (Note that when $\lambda=(0)$, $n_1=n_2$, (2) disappears.)
 \item For any $r$ such that $n_2+1 < r$, there exists $s=s(r) \geq 0$
       and $\mbf{i}=\mbf{i}(r)$ 
       such that $(\mbf{V}(r),x)$ is isomorphic to a direct sum of copies of 
       $(\Phi^+)^s I(\mbf{i})$, where
       $I(\mbf{i})$ is the injective representation corresponding to $\mbf{i}$; moreover
       if $n_2+1 <r<r'$, then either $s(r)>s(r')$ or $s(r)=s(r')$ and 
       $\mbf{i}(r) \neq \mbf{i}(r')$.
\end{enumerate}
Note that $\mbf{V}(r)=0$ for $|r|\gg 0$. Denote by $n(x)_0$ (resp. $n(x)_1$) 
the unique integer such that 
$\mbf{V}(n(x)_0)\neq 0$ and $\mbf{V}(n(x)_0-1)= 0$ 
(resp. $\mbf{V}(n(x)_1) \neq 0$ and $\mbf{V}(n(x)_1+1)=0$). 
For any $n(x)_0\leq r\leq n(x)_1$, define a function $\sigma_r: \mbf{Ind}\to \mathbb{N}$ 
as follows. 
\begin{itemize}
\item If $r\neq n_2+1$, $\sigma_r([P])=\sigma([P])$ 
      if $[P]\in\mbf{Ind}$ is the isoclass of  a direct summand of 
      $(\mbf{V}(r),x)$, otherwise $\sigma_r([P])=0$.
\item If $r=n_2+1$, $\sigma_r\equiv 0$.
\end{itemize}
Let $\lambda^r=(0)$ if $r\neq n_2+1$; $\lambda$ if $r=n_2+1$.
The pair $(\sigma_r,\lambda^r)$ is then in
$\varphi(\mbf{V}(r),\Omega)$, where $n(x)_0\leq r\leq n(x)_1$. 
Denote by $\mbf{X}(\sigma_r)$ 
the subvariety $\mbf{X}(\sigma_r,\lambda^r)$ 
in $\mbf{E}_{\mbf{V}(r),\Omega}$, where $n(x)_0\leq r\leq n(x)_1$.
Given any $\mbf{W}$ with the same dimension vector as $\mbf{V}(r)$,
by abuse of notation, we still denote by $\mbf{X}(\sigma_r)$ the
subvariety $\mbf{X}(\sigma_r,\lambda^r)$ in $\mbf{E}_{\mbf{W},\Omega}$.
Note that by definitions,
if $r\neq n_2+1$, $\mbf{X}(\sigma_r)$ is a
$\mbf{G}_{\mbf{V}(r)}$-orbit of any element $x$ in $\mbf{X}(\sigma_r)$. 
Furthermore, we have
\begin{lem}
$\mbf{X}(\sigma_r)$ is open in $\mbf{E}_{\mbf{V}(r),\Omega}$, if
$r\leq n_1$ or $\;n_2+1<r$.
\end{lem}
\begin{proof} 
Fix an element $x$ in $\mbf{X}(\sigma_r)$.
Since 
$\mrm{dim}\; O_x=\mrm{dim}\;\mbf{G}_{\mbf{V}(r)}-
\mrm{dim}\;\mrm{St}_{\mbf{G}_{\mbf{V}(r)}}(x)$,
where $\mrm{St}_{\mbf{G}_{\mbf{V}(r)}}(x)$ is the stabilizer of $x$ in
$\mbf{G}_{\mbf{V}(r)}$. We have
\[
\mrm{dim}\;\mbf{E}_{\mbf{V}(r),\Omega}-\mrm{dim}\;O_x=\mrm{dim}\;\mbf{E}_{\mbf{V}(r),\Omega}-
\mrm{dim}\;\mbf{G}_{\mbf{V}(r)}+\mrm{dim}\;\mrm{St}_{\mbf{G}_{\mbf{V}(r)}}(x).
\]
Every element in $\mrm{St}_{\mbf{G}_{\mbf{V}(r)}}(x)$ 
can be regarded as an automorphism of $(\mbf{V}(r),x)$, 
so $\mrm{St}_{\mbf{G}_{\mbf{V}(r)}}(x)=\mrm{Aut}((\mbf{V}(r),x))$. 
Thus, $\mrm{St}_{\mbf{G}_{\mbf{V}(r)}}(x)$ has the same dimension as 
the Lie algebra of $\mrm{Aut}(\mbf{V}(r),x)$. 
So 
\[
\mrm{dim}\;\mbf{E}_{\mbf{V}(r),\Omega}-\mrm{dim}\;O_x
=\mrm{dim}\;\mbf{E}_{\mbf{V}(r),\Omega}-
\mrm{dim}\;\mbf{G}_{\mbf{V}(r)}+\mrm{dim}\;\mrm{Hom}((\mbf{V}(r),x),(\mbf{V}(r),x)).
\]
Now by 2.4 (1), we have 
\[
\mrm{dim}\;\mbf{E}_{\mbf{V}(r),\Omega}-\mrm{dim}\;O_x=
\mrm{dim}\;\mrm{Ext}^1((\mbf{V}(r),x),(\mbf{V}(r),x)).
\]
But $\mrm{Ext}^1((\mbf{V}(r),x),(\mbf{V}(r),x))=0$ by Lemma 2.6 A. We have
$\mrm{dim}\;\mbf{E}_{\mbf{V}(r),\Omega}=\mrm{dim}\;O_x$. So Lemma 5.4 holds.
\end{proof}

\subsection{}
Fix an element of $x$ in $\mbf{X}(\sigma,\lambda)$ as in 5.3.
Let 
$\boldsymbol{\nu}(x)=(\underline{\mrm{dim}}\;\mbf{V}(n(x)_0),
\cdots,\mrm{dim}\;\mbf{V}(n(x)_0+r),\cdots, \underline{\mrm{dim}}\;\mbf{V}(n(x)_1))$.
Let $\mathcal{F}_{\boldsymbol{\nu}(x)}$ be the variety of all flags of type 
$\boldsymbol{\nu}(x)$. 
Let
$\tilde{\mathcal{F}}_{\boldsymbol{\nu}(x)}$ be the 
variety consisting of all pairs $(x,\mbf{f})$, where $x\in \mbf{E}_{\mbf{V},\Omega}$ and
$\mbf{f}\in \mathcal{F}_{\boldsymbol{\nu}(x)}$, such that $\mbf{f}$ is $x$-stable.
Let $\tilde{\mathcal{F}}'$ (resp. $\tilde{\mathcal{F}}''$) 
be the subvariety of $\tilde{\mathcal{F}}_{\boldsymbol{\nu}(x)}$ 
consisting of all pairs $(x,\mbf{f})$, 
(where $\mbf{f}=(\mbf{V}=\mbf{V}^0\supseteq \cdots\supseteq\mbf{V}^N=0)$, 
$N=n(x)_1-n(x)_0+1$) 
such that the induced representation 
$(\mbf{V}^r/\mbf{V}^{r+1}, x)$ is in $\mbf{X}(\sigma_{n(x)_0+r})$ 
(resp. $\overline{\mbf{X}(\sigma_{n(x)_0+r}))}$), for any $r$. Note that 
$\mbf{X}(\sigma_{n(x)_0+r})\subseteq\mbf{E}_{\mbf{V}^r/\mbf{V}^{r+1},\Omega}$ is well-defined,
by definitions.
Since $\mbf{X}(\sigma_n)$ is open in its closure for any $n$,
$\tilde{\mathcal{F}}'$ is an open subvariety of $\tilde{\mathcal{F}}''$. Let 
$\boldsymbol{\pi'}:\tilde{\mathcal{F}}' \rightarrow \mbf{E}_{\mbf{V},\Omega}$ and 
$\boldsymbol{\pi''}:\tilde{\mathcal{F}}'' \rightarrow \mbf{E}_{\mbf{V},\Omega}$ 
be the first projections. $\boldsymbol{\pi}''$ is a proper morphism. We have 
\begin{lem}
 \begin{enumerate}
 \item The restriction of $\boldsymbol{\pi'}$ on 
       $\boldsymbol{\pi'}^{-1}(\mbf{X}(\sigma,\lambda))$,
       $\boldsymbol{\pi'}:\boldsymbol{\pi'}^{-1}(\mbf{X}(\sigma,\lambda))\to 
       \mbf{X}(\sigma,\lambda)$, is an isomorphism;
 \item $\tilde{\mathcal{F}}'$ and 
       $\tilde{\mathcal{F}}''$ are irreducible; they have the same dimension equal to
       $\mrm{dim}\;\mbf{X}(\sigma,\lambda)$;
 \item $\boldsymbol{\pi}''(\tilde{\mathcal{F}}'')=\overline{\mbf{X}(\sigma,\lambda)}$.
 \end{enumerate}
\end{lem}
\begin{proof}
First, write $\mbf{V}=\oplus_{n\in \mathbb{Z}}\mbf{V}(n)$ and order
$\mbf{V}(r)$ as in Section 5.3. 
Set $\mbf{V}[r]=\oplus_{m:m\geq  r}\mbf{V}(m)$ for any $r$.
The flag 
$\mbf{f}_x=(\mbf{V}[n(x)_0]\supseteq \cdots
\supseteq\mbf{V}[n(x)_0+r]\supseteq\cdots \supseteq \mbf{V}[n(x)_1]\supseteq 0)$ is 
then in $\mathcal{F}_{\boldsymbol{\nu}(x)}$. 
So $(x,\mbf{f}_x)$ is in $\boldsymbol{\pi}'^{-1}(x)$, by definitions. 
To prove (1), it suffices to prove that $\mbf{f}=\mbf{f}_x$, 
for any $x$-stable flag 
$\mbf{f}=(\mbf{V}=\mbf{V}^0\supseteq\mbf{V}^1\supseteq\cdots\supseteq \mbf{V}^N=0)$ in 
$\mathcal{F}_{\boldsymbol{\nu}(x)}$.
Now consider the following short exact sequence of representations of $\mbf{Q}$:
\[(\star) \hspace{1cm} 
0\to (\mbf{V}^{r+1},x)\to (\mbf{V}^r,x) \to (\mbf{V}^r/\mbf{V}^{r+1},x)\to 0.\]
By assumption, we have 
$(\mbf{V}^r/\mbf{V}^{r+1},x)\simeq (\mbf{V}(n(x)_0+r),x)$ for all $r \neq n_2-n(x)_0+1$ and 
$(\mbf{V}^{n_2-n(x)_0+1}/\mbf{V}^{n_2-n(x)_0+2},x)$ is a regular representation whose 
indecomposable summands are in different tubes,
but not in the tubes where $(\mbf{V}^r/\mbf{V}^{r+1},x)$ is in, 
for $n_1-n(x)_0<r \leq n_2-n(x)_0$. 
So we have
\[
\mathrm{Ext}^1((\mbf{V}^{r}/\mbf{V}^{r+1},x),(\mbf{V}^{r'}/\mbf{V}^{r'+1},x))=0
\hspace{1cm} \mbox{for any} \;\;r < r',
\]
by Lemma 2.6 A. By induction, ($\star$) is split for any $r$. 
So for  each $r$, there exists $x$-stable $\mbf{I}$-graded $k$-vector space $\mbf{W}^r$ 
such that
$\mbf{V}^r = \mbf{V}^{r+1}\oplus \mbf{W}^r$. We then have $\mbf{V}=\oplus_r\mbf{W}^r$ and 
$(\mbf{W}^r,x) \simeq (\mbf{V}(n(x)_0+r),x)$ 
for all $r \neq n_2-n(x)_0+1$. 
Thus we have automatically 
$(\mbf{W}^r,x) \simeq (\mbf{V}(n(x)_0+r),x)$, for $r=n_2-n(x)_0+1$. 
Hence we can find an automorphism $\phi$ of $(\mbf{V},x)$ 
which maps $\mbf{V}(n(x)_0+r)$ onto $\mbf{W}^r$, for all $r$. Let
$\phi_{r,r'}:\mbf{V}(r) \rightarrow \mbf{V}(r')$ 
be the composition of the $\mbf{I}$-graded linear maps 
$\mbf{V}(r)\to \mbf{V} \overset{\phi}\to \mbf{V} \to \mbf{V}(r')$. 
Since $\phi$ is compatible with $x$, 
it defines a morphism in $Rep(\mbf{Q})$, 
i.e. $\phi_{r,r'} \in \mrm{Hom}((\mbf{V}(r),x),(\mbf{V}(r'),x))$. 
So by Lemma 2.6 B, it must be zero, whenever $r>r'$. Thus we have 
\[
\phi(\mbf{V}(r)) \subseteq \bigoplus_{r':r' \geq r }\mbf{V}(r')=\mbf{V}[r]
\hspace{.1cm}
.
\]
Hence $\phi$ maps $\mbf{V}[r]$ into itself, for any $r$. Since $\phi$ is an isomorphism, 
we have  
$\phi(\mbf{V}[r])=\mbf{V}[r]$. 
On the other hand, from the definition of $\phi$, 
we have $\phi(\mbf{V}[n(x)_0+r])=\mbf{V}^r$. 
Thus $\mbf{V}^r=\mbf{V}[n(x)_0+r]$ for all $r$. (1) is proved.

Second, we consider the second projection 
$p_2':\tilde{\mathcal{F}}' \rightarrow \mathcal{F}_{\boldsymbol{\nu}(x)}$
(resp. 
$p_2'':\tilde{\mathcal{F}}'' \rightarrow \mathcal{F}_{\boldsymbol{\nu}(x)}$). 
This map is $\mbf{G_V}$-equivariant and $\mbf{G_V}$ acts transitively on 
$\mathcal{F}_{\boldsymbol{\nu}(x)}$. To prove (2), 
it's enough to prove that $(p_2')^{-1}(\mbf{f})$ and 
$(p_2'')^{-1}(\mbf{f})$ are irreducible varieties of dimension equal to 
$\mrm{dim}\;\mbf{X}(\sigma,\lambda)-\mrm{dim}\;\mathcal{F}_{\boldsymbol{\nu}(x)}$,
for any flag
$\mbf{f}=(\mbf{V}=\mbf{V}^0 \supseteq \cdots \supseteq \mbf{V}^N=0)$ in
$\mathcal{F}_{\boldsymbol{\nu}(x)}$. 
Now fix $\mbf{f}$ to be 
$\mbf{f}_x=
(\mbf{V}=\mbf{V}[n(x)_0] \supseteq \cdots \supseteq \mbf{V}[n(x)_1] \supseteq 0)$.  
We have 
\[
(p_2')^{-1}(\mbf{f}) \simeq
\prod_{r}\mbf{X}(\sigma_r)
\times
\prod_{r<r',h \in \Omega}\mrm{Hom}(\mbf{V}(r)_{h'},\mbf{V}(r')_{h''}),
\] 
and 
\[
(p_2'')^{-1}(\mbf{f})\simeq
\prod_r \overline{\mbf{X}(\sigma_r)}
\times
\prod_{r<r',h \in \Omega}\mrm{Hom}(\mbf{V}(r)_{h'},\mbf{V}(r')_{h''}).
\] 
So we have $(p_2')^{-1}(\mbf{f})$ and $(p_2'')^{-1}(\mbf{f})$ 
are irreducible of the same dimension and  
\begin{itemize}
\item[1.] $\mrm{dim}\;(p_2')^{-1}(\mbf{f})=\sum_{r}\mrm{dim}\;\mbf{X}(\sigma_r)+
          \sum_{r<r',\; h \in \Omega}\mrm{dim}\;\mrm{Hom}(\mbf{V}(r)_{h'},\mbf{V}(r')_{h''}).$
\end{itemize}
By Proposition 5.2 (1), we have 
\begin{itemize}
\item[2.] 
$\mrm{dim}\;\mbf{X}(\sigma,\lambda)=q+\mrm{dim}\;\mbf{G_V}-\mrm{dim}\;
\mrm{St}_{\mbf{G_V}}(x).$
\end{itemize}
Since $\mrm{St}_{\mbf{G_V}}(x)$ has the same 
dimension as the Lie algebra of $\mrm{Aut}(\mbf{V},x)$;
so
\[
\mrm{dim}\;\mrm{St}_{\mbf{G_V}}(x)=
\sum_{r,r'}\mrm{dim}\;\mrm{Hom}((\mbf{V}(r),x),(\mbf{V}(r'),x)).
\]
By Lemma 2.6 B, we have 
\[
\mrm{Hom}((\mbf{V}(r),x),(\mbf{V}(r'),x))=0 \hspace{2cm}\mbox{if} \; \;r>r'.
\]
Hence,  
\begin{itemize}
\item[3.]
         $\mrm{dim}\;\mrm{St}_{\mbf{G_V}}(x)=\sum_{r \leq r'}\mrm{dim}\;
         \mrm{Hom}((\mbf{V}(r),x),(\mbf{V}(r'),x)) \;.$
\end{itemize}
On the other hand, since $\mbf{G_V}$ acts on 
$\mathcal{F}_{\boldsymbol{\nu}(x)}$ transitively, 
$\mrm{dim}\;\mathcal{F}_{\boldsymbol{\nu}(x)}=
\mrm{dim}\;\mbf{G_V}-\mrm{dim}\;\mrm{St}_{\mbf{G_V}}(\mbf{f})$, 
where $\mrm{St}_{\mbf{G_V}}(\mbf{f})$ is the stabilizer of $\mbf{f}$ in $\mbf{G_V}$. Thus
\begin{itemize}
\item[4.]
$\mrm{dim}\;\mathcal{F}_{\boldsymbol{\nu}(x)}=\mrm{dim}\;\mbf{G_V}-
\sum_{r \leq r',\;\mbf{i \in I}}\mrm{dim}\;
\mrm{Hom}(\mbf{V}(r)_{\mbf{i}},\mbf{V}(r')_{\mbf{i}}).$
\end{itemize}
From (1), (2), (3) and (4) above, we see that 
\[\hspace{-9cm}
\mrm{dim}\;(p_2')^{-1}(\mbf{f})-\mrm{dim}\;\mbf{X}(\sigma,\lambda)+
\mrm{dim}\;\mathcal{F}_{\boldsymbol{\nu}(x)}
\]
\[
\hspace{-4cm}
=
\sum_r\mrm{dim}\;\mbf{X}(\sigma_r)+
\sum_{r<r',h \in
  \Omega}\mrm{dim}\;\mrm{Hom}(\mbf{V}(r)_{h'},\mbf{V}(r')_{h''})-q
-\mrm{dim}\;\mbf{G_V}
\]
\[
\hspace{-1cm}
+\sum_{r \leq r'}\mrm{dim}\;\mrm{Hom}((\mbf{V}(r),x),(\mbf{V}(r'),x))+
\mrm{dim}\;\mbf{G_V}-\sum_{r \leq r',\;\mbf{i \in I}}
\mrm{dim}\;\mrm{Hom}(\mbf{V}(r)_{\mbf{i}},\mbf{V}(r')_{\mbf{i}}).
\]
It remains to show that the last expression is zero. By 2.4 (1), 
the last expression is equal to
\[
\sum_r \mrm{dim}\;\mbf{X}(\sigma_r)-q+\sum_{r<r'}\mrm{dim} \;
\mrm{Ext}^1((\mbf{V}(r),x),(\mbf{V}(r'),x))
\]
\[\hspace{1cm}
+\sum_r \mrm{dim}\;\mrm{Hom}((\mbf{V}(r),x),(\mbf{V}(r),x))-
\sum_{\mbf{i},\; r}\mrm{dim}\;
\mrm{Hom}(\mbf{V}(r)_{\mbf{i}},\mbf{V}(r)_{\mbf{i}}).
\] 
Thus, by Lemma 2.6 A, this is equal to 
\[
\sum_r\mrm{dim}\;\mbf{X}(\sigma_r)-q+\sum_r\mrm{dim}\;\mrm{Hom}((\mbf{V}(r),x),
(\mbf{V}(r),x))-\sum_r\mrm{dim}\;\mbf{G}_{\mbf{V}(r)}.
\]
To prove that this is zero, it's enough to show that 
\[
\mrm{dim}\;\mbf{X}(\sigma_r)=\mrm{dim}\;\mbf{G}_{\mbf{V}(r)}-
\mrm{dim\; Hom}((\mbf{V}(r),x),(\mbf{V}(r),x)),
\]
for all $r\neq n_2-n(x)_0+1$ and 
\[
\mrm{dim}\;\mbf{X}(\sigma_r)=\mrm{dim}\;\mbf{G}_{\mbf{V}(r)}-
\mrm{dim\; Hom}((\mbf{V}(r),x),(\mbf{V}(r),x)) +q,
\]
for $r=n_2-n(x)_0+1$. 
These follows immediately from the definitions and from 
Proposition 5.2 (1). This completes the proof of (2).

Finally, since $\boldsymbol{\pi}''$ is proper, 
we see from (2) that the image of $\boldsymbol{\pi}''$ 
is a closed irreducible subset of $\mbf{E}_{\mbf{V},\Omega}$ of dimension
$\leq \mrm{dim}\;\mbf{X}(\sigma,\lambda)$. 
This image contains $\mbf{X}(\sigma,\lambda)$ and therefore 
it has dimension equal to $\mrm{dim}\;\mbf{X}(\sigma,\lambda)$. 
(3) holds.
\end{proof}

\subsection{}
Consider the cyclic quiver $\mbf{Q}_T$ corresponding to a tube $T$ 
(of period $p(T)$), as defined in Proposition 3.4. 
Given any pair $(r,m)$, 
where $r\in \mathbb{Z}/p\mathbb{Z}$ and $m\in \mathbb{N}$, 
there is a unique (up to isomorphism) indecomposable representation 
of $\mbf{Q}_T$ as follows.
\[V_{r,m}= \cdots 0 \rightarrow k \rightarrow \cdots \rightarrow k \rightarrow 0 
\cdots,\]
where the sequence starts at $r$ (mod $p(T)$) and the length is $m$.  
The set $\{V_{r,m}|\;r\in \mathbb{Z}/p(T)\mathbb{Z},m\in \mathbb{N}\}$ 
forms a complete list of indecomposable objects of $Nil(\mbf{Q}_{T})$ 
(up to isomorphism). 
We denote by $\mbf{V}_{T,r,m}$ the indecomposable object of $Rep(T)$ 
corresponding to $V_{r,m}$ under the functors $F$ and $G$ in the proof of Proposition 3.4.
From the definitions, we have 
\[
\underline{\mrm{dim}}\;\mbf{V}_{T,r,m}=
\sum^{r+m-1}_{s=r} \underline{\mrm{dim}} \;(\Phi^+)^{s}R_T,
\]
where $R_T$ is the fix regular simple in $T$.
For each $\nu\in\mathbb{N}[\mbf{I}]$, 
we fix a representative, $\mbf{V}_{\nu}$ of the $\mbf{I}$-graded $k$-vector spaces 
of dimension vector $\nu$. 
Consider the formal power series
\[
\Pi=\sum_{\nu\in\mathbb{N}[\mbf{I}]}|\varphi(\mbf{V}_{\nu},\Omega)|\;
X^{{\mrm{dim}}\;\mbf{V}_{\nu}}.
\]
From the definitions and Theorem 2.7, we see that 
\[\Pi=\prod_{\alpha \in \mbf{R}_{+}-\mbf{R}^{0}_{+}}(1-X^{\mrm{dim}(\alpha)})^{-1} 
\prod_{T,m}(\prod_{r} (1-X^{\mrm{dim}\;\mbf{V}_{T,r,m}})^{-1} - 
\prod_{r} X^{\mrm{dim}\;\mbf{V}_{T,r,m}} 
\prod_{r}(1-X^{\mrm{dim}\;\mbf{V}_{T,r,m}})^{-1}) \]

\[\times \prod_{T,s}((1-X^{sN})^{-p(T)}-X^{p(T)sN}(1-X^{sN})^{-p(T)}) 
\sum_{t \geq 0} q(t)X^{tN},\]
where $R_+$ is the set of all positive real roots, 
$R_+^0$ is the subset of $R_+$ consisting of all positive real roots 
whose corresponding indecomposable representations are regulars. 
$T$ runs over the set, $\mbf{T}$, of all tubes of period $\neq 1$, 
$m$ runs over the integers $\geq 1$, not divisible by $p(T)$, 
$r$ runs over $\mathbb{Z}/p(T)\mathbb{Z}$, 
$s$ runs over the integers $\geq 1$, 
and $q(t)$ is the number of partitions of t, 
$N=|\mbf{I}|$. Notice that 
$\mrm{dim}\;\mbf{V}_{T,r,m}=
\sum^{r+m-1}_{u=r}\mrm{dim}\;(\Phi^+)^{u}R_T$. 
We have
\[
\Pi=
\prod_{\alpha \in \mbf{R}_{+}-\mbf{R}^{0}_{+}}(1-X^{\mrm{dim}\;\alpha})^{-1} 
\prod_{T,m}\prod_{r}(1-X^{\sum^{r+m-1}_{u=r}\mrm{dim}\;(\Phi^+)^{u}R})^{-1}
\prod_{T,m}(1-\prod_{r}X^{\sum^{r+m-1}_{u=r}\mrm{dim}\;(\Phi^+)^{u}R})
\]
\[
\hspace{-4.7cm}
\times
\prod_{T,s}(1-X^{sN})^{-p(t)}
\prod_{T,s}(1-X^{p(T)sN})
\sum_{t \geq 0}q(t)X^{tN},
\]
\[
\hspace{-3.8cm}
=
\prod_{\alpha \in \mbf{R}_{+}-\mbf{R}^0_{+}}(1-X^{\mrm{dim}\;\alpha})^{-1}
\prod_{\alpha \in \mbf{R}^0_+}(1-X^{\mrm{dim}\;\alpha})^{-1}
\prod_{T,m}(1-X^{mN})
\]
\[
\hspace{-3.8cm}
\times
\prod_{T,s}(1-X^{sN})^{-p(T)}
\prod_{T,s}(1-X^{p(T)sN})
\prod_{s}(1-X^{sN})^{-1},
\]

\[
=
\prod_{\alpha \in \mbf{R}_{+}}(1-X^{\mrm{dim}\;\alpha})^{-1}
\prod_{s}(1-X^{sN})^{|\mbf{T}|}
\prod_{s}(1-X^{sN})^{-\sum{p(T)}}
\prod_{s}(1-X^{sN})^{-1}.
\]
Since $\sum_{T \in \mbf{T}}(p(T)-1)=N-2$ ([DR]), it follows that 

\[
\hspace{-4.5cm}
(\mbox{a})
\hspace{2cm}
\Pi
=
\prod_{\alpha \in \mbf{R}_{+}}(1-X^{\mrm{dim}\;\alpha})^{-1}
\prod_{s}(1-X^{sN})^{-N+1}
\hspace{.1cm}.
\]
On the other hand, by the Poincar\'{e}-Birkhoff-Witt Theorem, 
we have the equality of formal power series
\[
\hspace{-1.4cm}
(\mbox{b})
\hspace{2cm}
\sum_{\nu\in\mathbb{N}[\mbf{I}]}\mrm{dim}\;U^{-}_{\nu}\;X^{\mrm{dim}\;\nu}
=
\prod_{\alpha \in \mbf{R}_{+}}(1-X^{\mrm{dim}\;\alpha})^{-1}
\prod_{s}(1-X^{sN})^{-N+1}.
\]
By comparing ($a$) with ($b)$, we have
\[
\hspace{-7cm}
(\mbox{c})
\hspace{1cm}
\sum_{\nu:\; \mrm{dim} \; \nu =d}
(\mrm{dim}\;U^-_{\nu} - |\varphi(\mbf{V_{\nu}}, \Omega)|)=0,
\]
for any $d \geq 0$.

\begin{thm} For any $\mbf{I}$-graded $k$-vector space $\mbf{V}$, 
the map $(\sigma,\lambda) \mapsto N(\sigma,\lambda)$ is a bijection 
$\varphi(\mbf{V},\Omega) \simeq \mbf{Irr}(\Lambda_{\mbf{V}})$.
\end{thm}
\begin{proof}
By Proposition 5.2 (3), the map is injective.
For any $\mbf{V}$ such that $\underline{\mrm{dim}}\;\mbf{V}=\nu$, 
by Proposition 5.2 (2),we have 
\[
|\varphi(\mbf{V},\Omega) \leq |\mbf{Irr}(\Lambda_{\mbf{V}})|.
\]
By Lemma 4.3, we have
\[|\mbf{Irr}(\Lambda_{\mbf{V}})| \leq \mrm{dim}\;
\mathit{F}_{\mbf{V}}.\]
By Lemma 4.5, we have
\[\mrm{dim}\;\mathit{F}_{\mbf{V}} = \mrm{dim}\;U^-_{\nu}.\]
Thus, 
$\mrm{dim} \;U^-_{\nu}-|\varphi(\mbf{V},\Omega)| \geq 0$. 
Combining this with 5.7 (c), we have 
$\mrm{dim} \;U^-_{\nu}=|\varphi(\mbf{V},\Omega)|$. 
This implies that all inequalities above are equalities. In particular,
$|\varphi(\mbf{V},\Omega)|=|\mbf{Irr}(\Lambda_{\mbf{V}})|$.
Since the two sets are finite with the same number of elements and the
map is injective, the Theorem holds.
\end{proof}
\section{Cyclic quivers}
In this section, 
let $\mbf{Q}$ be a cyclic quiver defined in Proposition 3.4, with $p=N$. 
Denote by $Ind$ the set of all isomorphic classes of 
indecomposable nilpotent representations of $\mbf{Q}$. 
$Ind$ consists of the objects $[V_{r,m}]$, 
where $r\in \mathbb{Z}/N\mathbb{Z}$ and 
$m$ is an integer $\geq 1$ (See 5.7). 
A representation $(V,x)$ of $\mbf{Q}$ is called $aperiodic$ if, 
for any $m$, 
not all the indecomposables $V_{1,m},\cdots, V_{N,m}$ are 
direct summands of $(V,x)$. 
Given $V$, a $\mathbb{Z}/N\mathbb{Z}$-graded $k$-vector space, 
we denote by $\varphi(V)$ the set of all functions 
$\sigma:Ind \rightarrow \mathbb{N}$ 
such that the following three conditions are satisfied.
\begin{enumerate}
 \item $\sigma$ has finite support;
 \item $\prod_{r\in\mathbb{Z}/N\mathbb{Z}} \sigma([V_{r,m}])=0$, for any $m\geq 1$;
 \item $\sum_{[P]\in Ind} \sigma([P]) \;\underline{\mrm{dim}}\;P=\underline{\mrm{dim}}\;V$.
\end{enumerate}
Given $\sigma \in \varphi(V)$, 
define $\mbf{X}(\sigma)$ to be the set of all $x \in \mbf{E}_{V,\Omega}$ 
such that $(V,x)$ is isomorphic to $\oplus_{[P]\in Ind}P^{\sigma([P])}$.
Note that $\mbf{X}(\sigma)$ is just 
the $\mbf{G}_V$-orbit of any element $x\in\mbf{X}(\sigma)$ and that 
$(V,x)$ is aperiodic, for any $x\in\mbf{X}(\sigma)$. 
Let $P_{\sigma}$ be the simple perverse sheaf on
$\mbf{E}_{V,\Omega}$ 
whose support is the closure of $\mbf{X}(\sigma)$
and whose restriction to $\mbf{X}(\sigma)$ 
is the $l$-adic constant sheaf $\bar{\mathbb{Q}}_l$, up to a shift. 
($l$ is a prime number invertible in $k$.)  
Let $P_V$ be the set of all isomorphic classes of 
simple perverse sheaves on $\mbf{E}_{V,\Omega}$ 
in the class defined in [L1, \S2], for a cyclic quiver. 
We have 
\begin{thm}
For any $V$, the map $\sigma \mapsto P_{\sigma}$ is a bijection $\varphi(V) \simeq P_V$.
\end{thm}
See [L2, 5.9] for a proof. 
Note that this Theorem says that the simple perverse sheaves in $P_V$ 
are $1-1$ corresponding to 
the orbits of the aperiodic representations in $\mbf{E}_{V,\Omega}$.

\section{None cyclic quivers of type $\tilde{\mbf{A}}_n$}
In this section, 
we assume that $\mbf{Q}=(\Gamma,\Omega)$ 
is an affine quiver of type $\tilde{A}_n$, but not a cyclic quiver. 
We also assume that $l$ is a prime number invertible in $k$. 
Given an algebraic variety $X$, 
denote by $\bar{\mathbb{Q}}_l$ 
the $l$-adic constant sheaf on this variety. 
Denote by $\mathcal{D}_c^b(X)$ the bounded derived category
of complexes of $l$-adic sheaves on $X$ over $k$. 
\subsection{}
Let $\mbf{V}$ be an $\mbf{I}$-graded vector space over $k$. 
Given 
$\boldsymbol{\nu}=(\nu^1, \cdots, \nu^m) \in \mathcal{Y}$ 
(see 4.2) such that 
$|\boldsymbol{\nu}|=\underline{\mrm{dim}}\;\mbf{V}$, 
define $\mathcal{F}_{\boldsymbol{\nu}}$ to be the variety of 
all flags of type $\boldsymbol{\nu}$. 
Denote by $\tilde{\mathcal{F}}_{\boldsymbol{\nu}}$ the 
variety consisting of all pairs $(x,\mbf{f})$ such that $x \in \mbf{E}_{\mbf{V},\Omega}$,  
$\mbf{f} \in \mathcal{F}_{\boldsymbol{\nu}}$ and $\mbf{f}$ is $x$-stable. 
Let 
$\boldsymbol{\pi_{\nu}}:\tilde{\mathcal{F}}_{\boldsymbol{\nu}} 
\rightarrow \mbf{E}_{\mbf{V},\Omega}$ 
be the first projection. 
It induces a right derived functor 
$(\boldsymbol{\pi_{\nu}})_!:\mathcal{D}_c^b(\tilde{\mathcal{F}}_{\boldsymbol{\nu}})\to
\mathcal{D}_c^b(\mbf{E}_{\mbf{V},\Omega})$.
Note that $\boldsymbol{\pi_{\nu}}$ is proper. 
So by the Decomposition Theorem in [BBD], 
$(\boldsymbol{\pi_{\nu}})_!(\bar{\mathbb{Q}}_l)$ is semisimple. 
Let $\mbf{P}_{\mbf{V},\Omega}$ be 
the set of isoclasses of simple perverse sheaves 
$\mathrm{P}$ on $\mbf{E}_{\mbf{V},\Omega}$ 
satisfying the following condition: 
$\mathrm{P}$ is a direct summand of 
$(\boldsymbol{\pi_{\nu}})_!(\bar{\mathbb{Q}}_l)$ 
for some $\boldsymbol{\nu} \in \mathcal{Y}$, up to a shift.
The main goal of this paper is 
to describe the elements in $\mbf{P}_{\mbf{V},\Omega}$ 
by specifying their supports and the corresponding local systems.

\subsection{}
In this subsection, we study some special cases.
First we have 
\begin{lem}
For any $\mbf{I}$-graded $k$-vector space $\mbf{V}$, 
let $d=\mrm{dim}(\mbf{E}_{\mbf{V},\Omega})$,
then the simple perverse sheaf $\bar{\mathbb{Q}}_l[d]$ 
on $\mbf{E}_{\mbf{V},\Omega}$ is in $\mbf{P}_{\mbf{V},\Omega}$.
\end{lem}
\begin{proof}
Given any quiver, as long as this quiver has no oriented cycles. 
We can arrange the vertices in $\mbf{I}$ in a sequence: 
$\mbf{i}_1,\cdots,\mbf{i}_N$ ($N=|\mbf{I}|$)  
such that $\mbf{i}_r$ is a source of the full subquiver $\mbf{Q}(r)$ 
with vertex set $\mbf{I}-\{\mbf{i}_1,\cdots,\mbf{i}_{r-1}\}$. 
For any $\mbf{V}$, let 
$\boldsymbol{\nu}=(\mrm{dim}\mbf{V}_{\mbf{i}_1} \;\mbf{i}_1 ,
\cdots,\mrm{dim}\mbf{V}_{i_N}\; \mbf{i}_N)$.
By definitions, 
$\mathcal{F}_{\boldsymbol{\nu}}$ consists of just a single flag. 
Moreover, this flag is stable under any $x$ in $\mbf{E}_{\mbf{V},\Omega}$. 
Thus,   
$\boldsymbol{\pi_{\nu}}:\tilde{\mathcal{F}}_{\boldsymbol{\nu}} \to \mbf{E}_{\mbf{V},\Omega}$ 
is an isomorphism. 
Therefore, we have 
$(\boldsymbol{\pi_{\nu}})_!(\bar{\mathbb{Q}}_l)[d]=
\bar{\mathbb{Q}}_l[d]\in \mbf{P}_{\mbf{V},\Omega}$.
\end{proof}

Second, we assume that $\mbf{V}$ is an $\mbf{I}$-graded $k$-vector space such that 
$\underline{\mrm{dim}}\;\mbf{V}=q \delta$ ($\delta=\sum_{\mbf{i\in I}}\mbf{i}$).
Define  $\mbf{X}(0)$ to be  the subvariety of
$\mbf{E}_{\mbf{V},\Omega}$ consisting of all $x$ such that 
\[
(\mbf{V},x) \simeq R_1 \oplus \cdots \oplus R_q,
\]
where $R_1,\cdots,R_q$ are nonisomorphic regular simples of period 1. 
Note that $\mbf{X}(0)$ is a special case of $\mbf{X}(\sigma, \lambda)$
with $\sigma\equiv 0$ and $\lambda=(1,\cdots,1)$ in Section 5.1. 
Moreover, the closure of $\mbf{X}(0)$ is $\mbf{E}_{\mbf{V}, \Omega}$. 
Denote by $\tilde{\mbf{X}}(0)$ the variety consists of 
all sequences $(x, [R_1],\cdots,[R_q])$
such that $(\mbf{V},x)\simeq R_1\oplus\cdots\oplus R_q$,
where $x\in \mbf{X}(0)$ and $[R_1],\cdots,[R_q]$ 
are isoclasses of nonisomorphic regular simples of period 1,.
The first projection $Pr_1:\tilde{\mbf{X}}(0)\to \mbf{X}(0)$ is an $S_q$-principle covering, 
where $S_q$ is the symmetric group  of $q$ letters. 
$S_q$ acts naturally on $(Pr_1)_{\star}(\bar{\mathbb{Q}}_l)$. 
Given an irreducible representation $\chi$ of $S_q$, 
denote by $L_{\chi}$ the direct summand of 
$(Pr_1)_{\star}(\bar{\mathbb{Q}}_l)$ 
on which $S_q$ acts the same way as the character of $\chi$.

Let $\mbf{P}_{0,\chi}$ be  the simple perverse sheaf on 
$\mbf{E}_{\mbf{V},\Omega}$ 
whose support is $\mbf{E}_{\mbf{V},\Omega}$
and whose restriction to $\mbf{X}(0)$ is $\mathit{L}_{\chi}$ 
(up to a shift). We then have:
\begin{lem}
$\mbf{P}_{0,\chi} \in \mbf{P}_{\mbf{V},\Omega}$
\end{lem}
\begin{proof}
Let $\mathcal{F}_{\boldsymbol{\delta}}$ be the variety consisting of 
all flags of type $\boldsymbol{\delta}$, 
where $\boldsymbol{\delta}=(\delta,\cdots,\delta)$
such that 
$|\boldsymbol{\delta}|=q\;\delta$. 
Let 
$\tilde{\mathcal{F}}_{\boldsymbol{\delta}}$ 
be the variety of all pairs $(x,\mbf{f})$ such that 
$x \in \mbf{E}_{\mbf{V},\Omega}$, 
$\mbf{f} \in \mathcal{F}_{\boldsymbol{\delta}}$ and $\mbf{f}$ is $x$-stable. 
Let 
$\boldsymbol{\pi_{\delta}}:\tilde{\mathcal{F}}_{\boldsymbol{\delta}} 
\rightarrow \mbf{E}_{\mbf{V},\Omega}$ 
be the first projection. 
If $x \in \mbf{X}(0)$ and 
$\mbf{f} \in \mathcal{F}_{\boldsymbol{\delta}}$ is $x$-stable, 
then $(\mbf{V}^r,x)$ is regular for any $r$, by Lemma 2.5. 
Fix  $x\in \mbf{X}(0)$,  
there are exactly $q!$ flags in $\mathcal{F}_{\boldsymbol{\delta}}$ 
that are $x$-stable.
So the restriction of $\boldsymbol{\pi_{\delta}}$ 
on $\boldsymbol{\pi_{\delta}}^{-1}(\mbf{X}(0))$
defines a covering 
$\boldsymbol{\pi_{\delta}}^{-1}(\mbf{X}(0)) \rightarrow \mbf{X}(0)$.
It's isomorphic to $\tilde{\mbf{X}}(0) \rightarrow \mbf{X}(0)$, 
by the definition of $\mbf{X}(0)$. 
We see therefore that some shift of $\mbf{P}_{0,\chi}$ 
is a direct summand of $(\boldsymbol{\pi_{\delta}})_!(\bar{\mathbb{Q}}_l)$. 
On the other hand, by [L1],
$\boldsymbol{\pi_{\delta}}_!(\bar{\mathbb{Q}}_l)$ 
is an iterated $\star$-product of $q$ simple perverse sheaves 
of the form $\bar{\mathbb{Q}}_l$ on 
$\mbf{E}_{\mbf{V_{\delta}},\Omega}$ up to a shift 
(where $\mbf{V}_{\delta}$ is an $\mbf{I}$-graded vector space over $k$ such that 
$\underline{\mrm{dim}}\;\mbf{V}_{\delta}=\delta$). 
These are certainly in 
$\mbf{P}_{\mbf{V}_{\boldsymbol{\delta}},\Omega}$,
by Lemma 7.3. 
Since some shift of $\mbf{P}_{0,\chi}$ 
is a direct summand in such an iterated $\star$-product, 
it's contained in $\mbf{P}_{\mbf{V},\Omega}$ by [L1, 3.2, 3.4]. The Lemma is proved.
\end{proof}
Finally, 
let $T$ be a tube of period $p \neq 1$. 
Suppose that $(\mbf{V},x)\in Rep(T)$ is aperiodic.
(I.e. not all the indecomposable representations
$M,\Phi^+M,\cdots,(\Phi^+)^{p-1}M$ 
is a direct summand of $(\mbf{V},x)$, for any $M\in T$.)  
Let $\mbf{P}$ be the simple perverse sheaf, 
whose support is the closure of the $\mbf{G_V}$-orbit $O_x$ of $x$ 
in $\mbf{E}_{\mbf{V},\Omega}$, and 
whose restriction to $O_x$ is $\bar{\mathbb{Q}}_l$, up to a shift. 
Then, we have 
\begin{lem}
$\mbf{P} \in \mbf{P}_{\mbf{V},\Omega}$.
\end{lem}

\begin{proof}
The proof has five steps.

Step 1.
We have an isomorphism of varieties 
$F:\{id\}\oplus\mbf{E}_{\mbf{V},\Omega_2}\to
\mbf{E}_{F(\mbf{V}),\Omega_T}$ (3.5 (a)).  
Let $O_{F(x)}$ be the $\mbf{G}_{F(\mbf{V})}$-orbit of $F(x)$, 
let $P_{F(x)}$ be the simple perverse sheaf on 
$\mbf{E}_{F(\mbf{V}),\Omega_T}$ 
whose support is the closure of $O_{F(x)}$ and 
whose restriction to $O_{F(x)}$ is $\bar{\mathbb{Q}}_l$ (up to a shift). 

The aperiodicity of $(\mbf{V},x)$ implies 
the aperiodicity of $(F(\mbf{V}),F(x))$ in $Nil(\mbf{Q}_T)$. 
So, by Theorem 6.1, $P_{F(x)}\in P_{F(\mbf{V})}$. 
More precisely, 
let $B$ be the variety of all flags 
$(F(\mbf{V})=W^0 \supseteq W^1 \supseteq \cdots \supseteq W^N=0)$, 
where $W^m$ is a $\mathbb{Z}/p\mathbb{Z}$-graded subspace of
$F(\mbf{V})$ 
such that $\mrm{dim}\;W^m/W^{m+1}=1$
for $m=0,\cdots,N-1$. 
Let $\tilde{B}$ be the variety of all pairs $(y,\mbf{f})$ such that 
$y \in \mbf{E}_{F(\mbf{V}),\Omega_T}$ , $\mbf{f} \in B$ and $\mbf{f}$ is $y$-stable. 
Let $\pi_0:\tilde{B} \rightarrow \mbf{E}_{F(\mbf{V}),\Omega_T}$ 
be the first projection.
Then, by [L1, $\S$2],

(a) {\em
$P_{F(x)}$ is a direct summand of $(\pi_0)_!(\bar{\mathbb{Q}}_l)[d_1]$, for some $d_1$.
}

Step 2. For any $\mathbb{Z}/p\mathbb{Z}$-graded subspace $U$ of $F(\mbf{V})$,
we have  an $\mbf{I}$-graded subspace of $\mbf{V}$, 
$G(U)=\oplus_{\mbf{i\in I}} \mbf{V_i}$, 
where $\mbf{V_i}=U_r$ if $\mbf{i} \in supp((\Phi^+)^r R_T)$ for some $r$. 
Let $\mathcal{F}'$ be the variety of all flags 
$(\mbf{V}=\mbf{V}^0 \supseteq \mbf{V}^1 \supseteq \cdots \supseteq \mbf{V}^N=0)$ 
such that 
$\underline{\mrm{dim}}\mbf{V}^r/\mbf{V}^{r-1} 
=\underline{\mrm{dim}}((\Phi^+)^s R_T)$ 
for some $s=s(r)$ and 
$\mbf{V_i}^r=\mbf{V_j}^r$ if $\{\mbf{i,j}\} \subseteq supp((\Phi^+)^t R_T)$ for some $t$. 
Define a map $G: B\to \mathcal{F}'$ by 
\[
G((F(\mbf{V})=W^0 \supseteq W^1 \supseteq \cdots \supseteq W^N=0)):=
(\mbf{V}=\mbf{V}^0\supseteq G(W^1)\supseteq \cdots \supseteq G(W^N)=0).
\]  
By definitions, $G$ is an isomorphism of varieties. 

Let $\tilde{\mathcal{F}}'$ be the variety consisting of all pairs $(x,\mbf{f})$ such that 
$x$ is in $\{id\}\oplus\mbf{E}_{\mbf{V},\Omega_2}$ (3.5), 
$\mbf{f}$ is in $\mathcal{F}'$ and $\mbf{f}$ is $x$-stable. 
Let 
$\pi':\tilde{\mathcal{F}}' \to \{id\}\oplus\mbf{E}_{\mbf{V},\Omega_2}$ 
be the first projection. 
Let $O'$ be the $\mbf{H_V}$-orbit of $G(F(x))$ in $\{id\}\oplus\mbf{E}_{\mbf{V},\Omega_2}$, 
and 
let $\mbf{P}'$ be the simple perverse sheaf on 
$\{id\}\oplus\mbf{E}_{\mbf{V},\Omega_2}$ 
whose support is the closure of $O'$ 
in $\{id\}\oplus\mbf{E}_{\mbf{V},\Omega_2}$ and 
whose restriction to $O'$ is $\bar{\mathbb{Q}}_l$ (up to a shift). 
Then by (a), we have

(b) {\em
$\mbf{P}'$  is a direct summand of $(\pi')_!(\bar{\mathbb{Q}}_l)[d_2]$, for some $d_2$.
}

Step 3. Let $\mathcal{F}$ be the variety of all flags 
$(\mbf{V}=\mbf{V}^0 \supseteq \mbf{V}^1 \supseteq \cdots \supseteq \mbf{V}^N =0)$ 
such that 
$\underline{\mrm{dim}}\mbf{V}^r /\mbf{V}^{r+1}$ =  
$\underline{\mrm{dim}}(\Phi^+)^s R_T$ 
for some $s=s(r)$ ($r=1,\cdots,N-1$). 
Let $\tilde{\mathcal{F}}''$ be the variety of all pairs $(x,\mbf{f})$ 
such that $x$ is in $\mrm{Aut}_{\mbf{V},\Omega_1}\oplus\mbf{E}_{\mbf{V},\Omega_2}$ (3.5), 
$\mbf{f}$ is in $\mathcal{F}$ and $\mbf{f}$ is $x$-stable. 
Given $x\in \{id\}\oplus\mbf{E}_{\mbf{V},\Omega_2}$, 
if a flag in $\mathcal{F}$ is $x$-stable, 
it then follows from the definitions that 
such a flag must automatically be contained in $\mathcal{F}'$. 
We then see that 
$\tilde{\mathcal{F}}'' =\mbf{G_V} \times^{\mbf{H_V}} \tilde{\mathcal{F}}'$ 
in the same way as 
$\mrm{Aut}_{\mbf{V},\Omega_1}\oplus\mbf{E}_{\mbf{V},\Omega_2}
=\mbf{G_V}\times^{\mbf{H_V}}(\{id\}\oplus\mbf{E}_{\mbf{V},\Omega_2})$ 
and 
$O_x =\mbf{G_V}\times^{\mbf{H_V}}O'$ in Lemma 3.6.

Let 
$\pi'':\tilde{\mathcal{F}}'' \to 
\mrm{Aut}_{\mbf{V},\Omega_1}\oplus\mbf{E}_{\mbf{V},\Omega_2}$ 
be the first projection. 
Let $\mbf{P}''$ be the simple perverse sheaf on 
$\mrm{Aut}_{\mbf{V},\Omega_1}\oplus\mbf{E}_{\mbf{V},\Omega_2}$, 
whose support is the closure of $O_x$ in
$\mrm{Aut}_{\mbf{V},\Omega_1}\oplus\mbf{E}_{\mbf{V},\Omega_2}$ 
and 
whose restriction to $O_x$ is $\bar{\mathbb{Q}}_l$ (up to a shift). 
Then by (b), we have

(c) {\em
$\mbf{P}''$ is a direct summand of $(\pi'')_!(\bar{\mathbb{Q}}_l)[d_3]$, for some $d_3$.
}

Step 4. Let $\tilde{\mathcal{F}}$ be the variety as follows. 
\[\tilde{\mathcal{F}}=\{ (x,\mbf{f})| \;x \in \mbf{E}_{\mbf{V},\Omega},\;
\mbf{f} \in \mathcal{F},\;\mbf{f} \;is\; x\emph{-}stable \}.
\]  
Let 
$\pi:\tilde{\mathcal{F}} \rightarrow \mbf{E}_{\mbf{V},\Omega}$ 
be the first projection.
Let 
$i:\mrm{Aut}_{\mbf{V},\Omega_1}\oplus
\mbf{E}_{\mbf{V},\Omega_2}\to \mbf{E}_{\mbf{V},\Omega}$ 
be the inclusion. 
We have $\pi_!(\bar{\mathbb{Q}}_l)=i_!(\pi'')_!(\bar{\mathbb{Q}}_l)$.
The direct image functor preserves finite direct sums. We deduce from
(c) that

(d) {\em  
$\mbf{P}$ is a direct summand of $\pi_!(\bar{\mathbb{Q}}_l)[d_4]$, for some $d_4$.
}

Step 5.
Let $\mathcal{Z}$ be the set of all sequences 
$\boldsymbol{\nu}=(\nu^1,\cdots,\nu^N)$
such that 
$|\boldsymbol{\nu}|=\underline{\mrm{dim}}\;\mbf{V}$ and $\nu^r$ is the
dimension vector of $(\Phi^+)^sR_T$ for some $s=s(r)$.
By definition, 
$\tilde{\mathcal{F}}=\cup_{\boldsymbol{\nu}\in\mathcal{Z}}\tilde{\mathcal{F}_{\nu}}$.
So  
$\pi_!(\bar{\mathbb{Q}}_l)=
\oplus_{\boldsymbol{\nu}\in\mathcal{Z}}(\boldsymbol{\pi_{\nu}})_!
(\bar{\mathbb{Q}}_l)$ (See 7.1 for the notation).
By [L1], $\boldsymbol{\pi_{\nu}}(\bar{\mathbb{Q}}_l)$ can be
decomposed into a direct sum of shifts of simple perverse sheaves in
$\mbf{P}_{\mbf{V}_{\boldsymbol{\nu}},\Omega}$, 
where
$\underline{\mrm{dim}}\;\mbf{V}_{\boldsymbol{\nu}}=|\boldsymbol{\nu}|$.
So $\pi_!(\bar{\mathbb{Q}}_l)$ can be decomposed into a direct sum of
shifts of simple perverse sheaves in $\mbf{P}_{\mbf{V},\Omega}$. 
Therefore, $\mbf{P}$ is in $\mbf{P}_{\mbf{V},\Omega}$.
Lemma 7.3 is proved.
\end{proof}

\subsection{}
In this subsection, let $\mbf{V}$ be an $\mbf{I}$-graded $k$-vector
space. We can give a general description of 
the elements in $\mbf{P}_{\mbf{V},\Omega}$. 
Given a pair $(\sigma, \lambda)\in \varphi(\mbf{V},\Omega)$, 
let $q=\sum_r\lambda_r$. 
Let $\tilde{\lambda}=(1,\cdots,1)$, with $q$ $1$s. 
Define $\mbf{X}(\sigma,\tilde{\lambda})$ as in Section 5.1, 
i.e. it's a subvariety of $\mbf{E}_{\mbf{V},\Omega}$ 
consists of all elements $x$ such that 
$(\mbf{V},x)=\oplus_{[P]\in \mbf{Ind}}P^{\sigma([P])}\oplus
R_1\oplus\cdots\oplus R_q$, 
where $R_1,\cdots,R_q$ are regular simples of period 1.
Define $\tilde{\mbf{X}}(\sigma,\tilde{\lambda})$ 
to be the variety consisting of all sequences $(x,[R_1],\cdots,[R_q])$,
where $x \in \mbf{X}(\sigma,\tilde{\lambda})$ and
$\{[R_1],\cdots,[R_q]\}$ is the set of all isoclasses of 
indecomposable summands of $(\mbf{V},x)$ that are regular simples of period 1. 
By definition, $\{[R_1],\cdots,[R_q]\}$ is uniquely determined by $x$ up to order. 
The first projection 
$Pr_1:\tilde{\mbf{X}}(\sigma) \rightarrow \mbf{X}(\sigma)$ 
is a principle covering with group $S_q$. 
Note that if $\lambda=(0)$, 
$\mbf{X}(\sigma,\tilde{\lambda})$ is a $\mbf{G_V}$-orbit
of any elements $x$ in $\mbf{X}(\sigma,\tilde{\lambda})$   
and 
$Pr_1$ is the identity map.
Since $\lambda=(\lambda_1\geq\cdots\geq\lambda_p)$ 
gives an irreducible representation of $S_q$, 
denote by $L_{\chi(\lambda)}$ the corresponding local system 
on $\mbf{X}(\sigma,\tilde{\lambda})$, via this covering.

Let $\mbf{P}_{\sigma,\chi(\lambda)}$ be 
the simple perverse sheaf on $\mbf{E}_{\mbf{V},\Omega}$,
whose support is the closure of $\mbf{X}(\sigma,\tilde{\lambda})$ 
and 
whose restriction to $\mbf{X}(\sigma,\tilde{\lambda})$ is 
$\mathit{L}_{\chi(\lambda)}$ (up to a shift).
We have 
\begin{prop}
$\mbf{P}_{\sigma,\chi(\lambda)} \in \mbf{P}_{\mbf{V},\Omega}$.
\end{prop} 
\begin{proof}
For any integer $r$, 
let 
$\mbf{X}(\sigma_r)=\mbf{X}(\sigma, \tilde{\lambda}^r)
\subseteq \mbf{E}_{\mbf{V}(r),\Omega}$ 
be defined as in Section 5.3. 
(Note that we consider the pair $(\sigma, \tilde{\lambda})$, not $(\sigma,\lambda)$.)
Define a local system, $\mathit{L}_r$, on $\mbf{X}(\sigma_r)$ 
to be  
$\bar{\mathbb{Q}}_l$ if $r\neq n_2+1$ and  
$\mathit{L}_{\chi(\lambda)}$ if $r=n_2+1$. 
(Note that 
the local system $\mathit{L}_{\chi(\lambda)}$ on
$\mbf{X}(\sigma_{n_2+1})$
is defined in a way similar to 
the restriction $L_{\chi}$ of $\mbf{P}_{0,\chi}$ on $\mbf{X}(0)$ in
Lemma 7.4)
Let $\mbf{P}_r$ be 
the simple perverse sheaf on $\mbf{E}_{\mbf{V}(r),\Omega}$, 
whose support is in the closure of $\mbf{X}(\sigma_r)$ and 
whose restriction to $\mbf{X}(\sigma_r)$ is $\mathit{L}_r$ (up to a shift). Then
\[
\hspace{-8cm}
\mbox{(a)}
\hspace{5cm}
\mbf{P}_r \in \mbf{P}_{\mbf{V}(r),\Omega},
\]
by Lemma 7.3 and Lemma 5.4 if $r\leq n_1$ or $n_2+1<r$; by Lemma 7.4
if $r=n_2+1$ and by Lemma 7.5 if $n_1<r\leq n_2$. 
From [L1, 3.2, 3.5], it then follows that 
the iterated $\star$-product of the  simple perverse sheaves
\[
\hspace{-6.8cm}
\mbox{(b)}
\hspace{5cm}
\mbf{P}_{n(x)_0} \star \cdots \star \mbf{P}_{n(x)_1}
\]
on $\mbf{E}_{\mbf{V},\Omega}$ is a direct sum of shifts of 
simple perverse sheaves in $\mbf{P}_{\mbf{V},\Omega}$. 
Hence it's enough to show that some shift of 
the simple perverse sheaf $\mbf{P}_{\sigma,\chi}$ 
is a direct summand of  the complex (b) above.
From the proof of Lemma 5.6 (2), we have a natural projection 
$\alpha:(p_2')^{-1}(\mbf{f})\to \prod_r\mbf{X}(\sigma_r)$. 
Let 
$i:(p_2')^{-1}(\mbf{f})\to \tilde{\mathcal{F}}_{\boldsymbol{\nu}(x)}'$
be the natural inclusion. 
Then $L'=i_{\star}(\alpha^{\star}(\otimes_rL_r))$ is 
a $\mbf{G_V}$-equivariant local system on
$\tilde{\mathcal{F}}_{\boldsymbol{\nu}(x)}'$,
where $\otimes$ is the external tensor product. 
Let $\mbf{P}'$ be the simple perverse sheaf on 
$\tilde{\mathcal{F}}_{\boldsymbol{\nu}(x)}''$ 
whose support is $\tilde{\mathcal{F}}_{\boldsymbol{\nu}(x)}''$ 
and whose restriction to 
$\tilde{\mathcal{F}}_{\boldsymbol{\nu}(x)}'$ is $L'$, up to a shift. 
Then by definition of the complex (b), 
\[
\mbf{P}_{n(x)_0}\star\cdots\star\mbf{P}_{n(x)_1}\simeq(\boldsymbol{\pi}'')_!\mbf{P}'[d],
\]
for some $d$. By Lemma 5.6 (2), we have
\[
\mrm{dim}\;(\tilde{\mathcal{F}}_{\boldsymbol{\nu}(x)}''-
\tilde{\mathcal{F}}_{\boldsymbol{\nu}(x)}')
<
\mrm{dim}\;\mbf{X}(\sigma,\tilde{\lambda}).
\]
Hence, 
\[
\mrm{dim}\;\boldsymbol{\pi}''
(\tilde{\mathcal{F}}_{\boldsymbol{\nu}(x)}''-\tilde{\mathcal{F}}_{\boldsymbol{\nu}(x)}')
<
\mrm{dim}\;\mbf{X}(\sigma,\tilde{\lambda}).
\]
Thus, the set
\[
\mbf{X}'=\mbf{X}(\sigma,\tilde{\lambda})-(\mbf{X}(\sigma,\tilde{\lambda})\cap 
\boldsymbol{\pi}''(\tilde{\mathcal{F}}_{\boldsymbol{\nu}(x)}''-
\tilde{\mathcal{F}}_{\boldsymbol{\nu}(x)}'))
\]
is an open dense subset of $\mbf{X}(\sigma,\tilde{\lambda})$.

By Lemma 5.6 (1), the restriction of $\boldsymbol{\pi}''$ is an isomorphism 
$(\boldsymbol{\pi}'')^{-1}\mbf{X}' \simeq \mbf{X}'$. 
Under this isomorphism, the restriction of the local system
$\mathit{L}'$ 
to the subset $(\boldsymbol{\pi}'')^{-1}(\mbf{X}')$ of
$\tilde{\mathcal{F}}_{\boldsymbol{\nu}(x)}'$ 
corresponds to a local system on $\mbf{X}'$ which 
can be seen to be just the restriction of the local system
$\mathit{L}_{\chi(\lambda)}$ defining 
$\mbf{P}_{\sigma,\chi(\lambda)}$ on $\mbf{X}(\sigma,\tilde{\lambda})$.

Thus, the cohomology sheaves of $(\boldsymbol{\pi}'')_!\mbf{P}'$ 
restricted to $\mbf{X}'$ are equal to 
$\mathit{L}_{\chi(\lambda)}|_{\mbf{X}'}$ in one degree 
and zero in all other degrees. 
Let $\mbf{P}''$ be the simple perverse sheaves 
whose support is the closure of $\mbf{X}'$ and 
whose restriction to $\mbf{X}'$ 
is $\mathit{L}_{\chi(\lambda)}|_{\mathit{X}'}$, up to a shift. 
Since $(\boldsymbol{\pi}'')_!\mbf{P}'$ 
is known to be 
a direct sum of shifts of simple perverse sheaves and 
$\mbf{X}'$ is open dense in the support of
$(\boldsymbol{\pi}'')_!\mbf{P}'$ 
(see Lemma 5.6 (3)), 
it follows that some shift of 
$\mbf{P}''$ is a direct summand of $(\boldsymbol{\pi}'')_!\mbf{P}'$.

We have clearly 
$\mbf{P}''=\mbf{P}_{\sigma,\chi(\lambda)}$. 
Proposition 7.3 is proved.
\end{proof}
\begin{cor}
The assignment $(\sigma,\lambda)\mapsto \mbf{P}_{\sigma,\chi(\lambda)}$
defines an injective map 
$\varphi(\mbf{V},\Omega)\to \mbf{P}_{\mbf{V},\Omega}$. 
\end{cor}
It's also surjective, as is shown in the following Theorem 7.10.

\subsection{}
Let $\mbf{U}^-$ be the negative part of the quantized enveloping algebra 
attached to a symmetric generalized Cartan matrix $C$ of type $\tilde{A}_n$.
When $n=1$, $\mbf{U}^-$ is a $\mathbb{Q}(v)$-algebra, 
where $v$ is an indeterminate, with two generators
$\mathit{F}_{\mbf{i}}$ and $\mathit{F}_{\mbf{j}}$
and two relations:
\[\mathit{F}_{\mbf{i}}\;\mathit{F}^3_{\mbf{j}}-
(v^2+1+v^{-2})\;\mathit{F}_{\mbf{j}}\;\mathit{F}_{\mbf{i}}\;\mathit{F}^2_{\mbf{j}}+
(v^2+1+v^{-2})\;\mathit{F}^2_{\mbf{j}}\;\mathit{F}_{\mbf{i}}\;\mathit{F}_{\mbf{j}}-
\mathit{F}^3_{\mbf{j}}\;\mathit{F}_{\mbf{i}}=0,
\]
\[
\mathit{F}_{\mbf{j}}\;\mathit{F}^3_{\mbf{i}}-
(v^2+1+v^{-2})\;\mathit{F}_{\mbf{i}}\;\mathit{F}_{\mbf{j}}\;\mathit{F}^2_{\mbf{i}}+
(v^2+1+v^{-2})\;\mathit{F}^2_{\mbf{i}}\;\mathit{F}_{\mbf{j}}\;\mathit{F}_{\mbf{i}}-
\mathit{F}^3_{\mbf{i}}\;\mathit{F}_{\mbf{j}}=0.
\]
When $n\geq 2$, $\mbf{U}^-$ is a $\mathbb{Q}(v)$-algebra with
generators $\mathit{F}_{\mbf{i}}$, 
where $\mbf{i}\in\mbf{I}$ with $|\mbf{I}|=n+1$, and 
relations:
\[
\mathit{F}_{\mbf{i}}\;\mathit{F}^2_{\mbf{j}}-
(v+v^{-1})\;\mathit{F}_{\mbf{j}}\;\mathit{F}_{\mbf{i}}\;\mathit{F}_{\mbf{j}}+
\mathit{F}^2_{\mbf{j}}\;\mathit{F}_{\mbf{i}}=0,
\;\; 
\emph{if} \;c_{\mbf{i,j}}=-1;
\]
\[\mathit{F}_{\mbf{i}}\;\mathit{F}_{\mbf{j}}=
\mathit{F}_{\mbf{j}}\;\mathit{F}_{\mbf{i}},
\;\;
\emph{if} \;c_{\mbf{i,j}}=0.
\]

Let $\mbf{F_V}$ be the $\mathbb{Q}(\mathit{v})$-vector space 
spanned by the elements in $\mbf{P}_{\mbf{V},\Omega}$. 
For each $\nu \in \mathbb{N}[\mbf{I}]$, 
we fix an $\mbf{I}$-graded $k$-vector space 
$\mbf{V}_{\nu}$ of dimension vector $\nu$. 
Let 
$\mbf{F}=\oplus_{\nu\in\mathbb{N}[\mbf{I}]}\mbf{F}_{\mbf{V}_{\nu}}$.
It's independent of the choice of $\mbf{V}_{\nu}$, 
since $\mbf{F}_{\mbf{V}_{\nu}}$ and
$\mbf{F}_{\mbf{V'}_{\nu}}$ are isomorphic if 
$\mbf{V}_{\nu}\simeq \mbf{V'}_{\nu}$.
From [L1, \S 3, \S 10], we know that 
there is a $\mathbb{Q}(v)$-algebra structure on $\mbf{F}$ 
such that  the map
$(\mathit{F}^{r_1}_{\mbf{i}_1}/[r_1]^!)\cdots(\mathit{F}^{r_m}/[r_m]^!)
\mapsto
v^{-d(\boldsymbol{\nu})}(\boldsymbol{\pi_{\nu}})_!(\bar{\mathbb{Q}}_l)$,
where
$d(\boldsymbol{\nu})=\mrm{dim}\;\tilde{\mathcal{F}_{\boldsymbol{\nu}}}$,
$[r]^!=\prod_{s=1}^{r}(v^s-v^{-s})/(v-v^{-1})$ for any $r$ and 
$\boldsymbol{\nu}=(r_1\mbf{i}_1,\cdots,r_m\mbf{i}_m)$,
defines an algebra isomorphism:
\[
\hspace{-10cm}
\mbox{(a)}
\hspace{4cm}
\mbf{U}^- \to \mbf{F}.
\]
There is a natural grading 
$\mbf{U}^-=\oplus_{\nu\in \mathbb{N}[\mbf{I}]}\mbf{U}^-_{\nu}$ 
and the isomorphism (a) respects the gradings. Hence,
\[
\hspace{-8.5cm}
\mbox{(b)} 
\hspace{4cm}
\mrm{dim}\;\mbf{U}^-_{\nu} = \mrm{dim}\;\mbf{F}_{\mbf{V}}, 
\]
for any $\mbf{V}$ such that $\underline{\mrm{dim}}\;\mbf{V}=\nu$.
Furthermore, we have 
\begin{thm}
For any  $\mbf{I}$-graded $k$-vector space $\mbf{V}$, 
the map 
$\varphi(\mbf{V},\Omega) \to \mbf{P}_{\mbf{V},\Omega}$
in Corollary 7.8 is a bijection.
\end{thm}
\begin{proof}
For any $\mbf{V}$, by Corollary 7.8,
\[
|\mbf{P}_{\mbf{V},\Omega}| \geq |\varphi(\mbf{V},\Omega)|.
\]
By 7.9 (b), we know that
\[
\mrm{dim}\;\mbf{U}_{\nu}^- = \mrm{dim}\;\mbf{F}_{\mbf{V}}=|\mbf{P}_{\mbf{V},\Omega}|,
\]
where $\underline{\mrm{dim}}(\mbf{V})=\nu$.
Since the $\mathbb{Q}$-algebra $U^-$ is 
a specialization of the $\mathbb{Q}(v)$-algebra
$\mbf{U}^-$, we have
\[
\mrm{dim}\;U_{\nu}^- \geq \mrm{dim}\;\mbf{U}_{\nu}^-.
\]
As in the proof of Theorem 5.8, we have
\[
|\varphi(\mbf{V},\Omega)|=\mrm{dim}\;U_{\nu}^-.\]
Combining all these together we have the Theorem.
\end{proof}
Remark 1. Theorem 7.10 gives a description of the canonical basis
elements in the quantized enveloping algebra of type $\tilde{A}_n$.

Remark 2. Theorem 7.10 is a generalization of Theorem 6.16 (b) in [L2] in
the case of type $\tilde{A}_n$. But Theorem 6.16 (b) in [L2] applies to
affine quivers of type $\tilde{D}_n$ and $\tilde{E}_n$ too, with some
restrictions on the orientations of the quivers. 
It will be interesting to generalize the description of the canonical
bases elements in [L2, Theorem 6.16] for arbitrary orientations in the
cases of type $\tilde{D}_n$ and $\tilde{E}_n$.


\begin{thebibliography}{99999}\frenchspacing
\bibitem[\sf B] {B} A.Borel {\em Linear algebraic groups}, GTM 126
\bibitem[\sf BB] {BB} S.Benner and M.C.R.Butler, 
        {\em The equivalence of certain functors occurring in the 
        representation theory of artin algebras and species}, 
        J. London Math. Soc., {\bf 14} (1976), 183-187.
\bibitem[\sf BBD]{BBD} A.A.Beilinson,J.Bernstein,and P.Deligne,
        {\em Faisceaux pervers},
        Ast\'{e}risque {\bf 100} (1982).
\bibitem[\sf BGP]{BGP} I.N.Bernstein, I.M.Gelfand and B.A.Ponomarev, 
        {\em Coxeter functor and Gabriel's theorem}, 
        Russian Math. Surveys, {\bf 28} (1976), 17-32.
\bibitem[\sf CB]{Crawley-Boevey} W.Crawley-Boevey,
        {\em Lectures on Representations of Quivers}, preprint.
\bibitem[\sf DR]{Dlab&Ringel} V.Dlab and C.M.Ringel,
        {\em Indecomposable representations of graphs and algebras}, 
        Memoirs Amer. Math. Soc.,{\bf 173} (1976), 1-57.
\bibitem[\sf L1]{lusztig1} G.Lusztig,
        {\em Quivers,perverse sheaves and quantized enveloping algebras}, 
        Jour. Amer. Math. Soc. {\bf 4} (1991), 365-421.
\bibitem[\sf L2]{lusztig2} G.Lusztig,
        {\em Affine quivers and canonical bases}, 
        Publ. Math. IHES {\bf 76} (1992), 111-163.
\bibitem[\sf L3]{lusztig3} G.Lusztig,
        {\em Introduction to Quantum Groups}, 
        Birkh\"{a}user in Math. {1993}.
\bibitem[\sf L4]{lusztig4} G.Lusztig,
        {\em Semicanonical bases arising from enveloping algebras}, 
        Adv. Math. {\bf 151} (2000), 129-139.
\bibitem[\sf R1]{Ringel1} C.M.Ringel,
        {\em Representations of $\mbf{K}$-species and bimodules}, 
        J. of Alg. {\bf 41} (1976), 269-302.
\bibitem[\sf R2]{Ringel2} C.M.Ringel,
        {\em The preprojective algebra of a tame quiver:The 
        irreducible components of the module varieties}, 
        Contemporary Mathematics, {\bf 229}, 1998.
\end{thebibliography}
\end{document}